\documentclass[12pt,fleqn]{article}

\usepackage[english]{babel}

\usepackage{amsmath,amssymb,amsthm,mathrsfs,amsfonts,dsfont}
\usepackage{indentfirst}

\usepackage[dvips]{graphicx}

\title{On the common points of two families of $N$-spheres in the flat $N+1$ dimensional space, each of which passes through the vertexes of a given $N$-simplex}
\author{Vassil K. Tinchev\thanks{E-mail: tintschev@phys.uni-sofia.bg}\\
{\footnotesize  Department of Theoretical Physics,
                Faculty of Physics, Sofia University,}\\
{\footnotesize  5 James Bourchier Boulevard, Sofia~1164, Bulgaria }\\}
\date{08.09.2014 г.}

\begin{document}

\maketitle

\textit{In memory of Mr Vlaikov---my math teacher in the fourth grade. Thanks to him I discovered the beauty of mathematics. I also dedicate this work to my parents, Ionka and Kiril, who gave me a happy childhood.}

\begin{abstract}
Let two distinct $N$-simplexes be given in an Euclidean or pseudo-Euclidean $N+1$ dimensional space as each is defined by the coordinates of its $N+1$ vertexes. We consider the two families of $N$-spheres passing through the vertexes of the given $N$-simplexes and the set of couples of $N$-spheres (one belonging to first family and the other to the second one). The elements of this set have at least one common point; moreover, it is such that for the angle $\alpha$ between the segments connecting that point and the centers of the corresponding $N$-spheres, there holds $\cos^{2}\alpha=const$ for each of the elements of the defined set of $N$-spheres. In the present work we find the geometric place of all these common points, including the special cases when $\cos^{2}\alpha$ is equal to $0$ or $1$.
\end{abstract}

\begin{figure}[h]
    \includegraphics[width=0.98\textwidth]{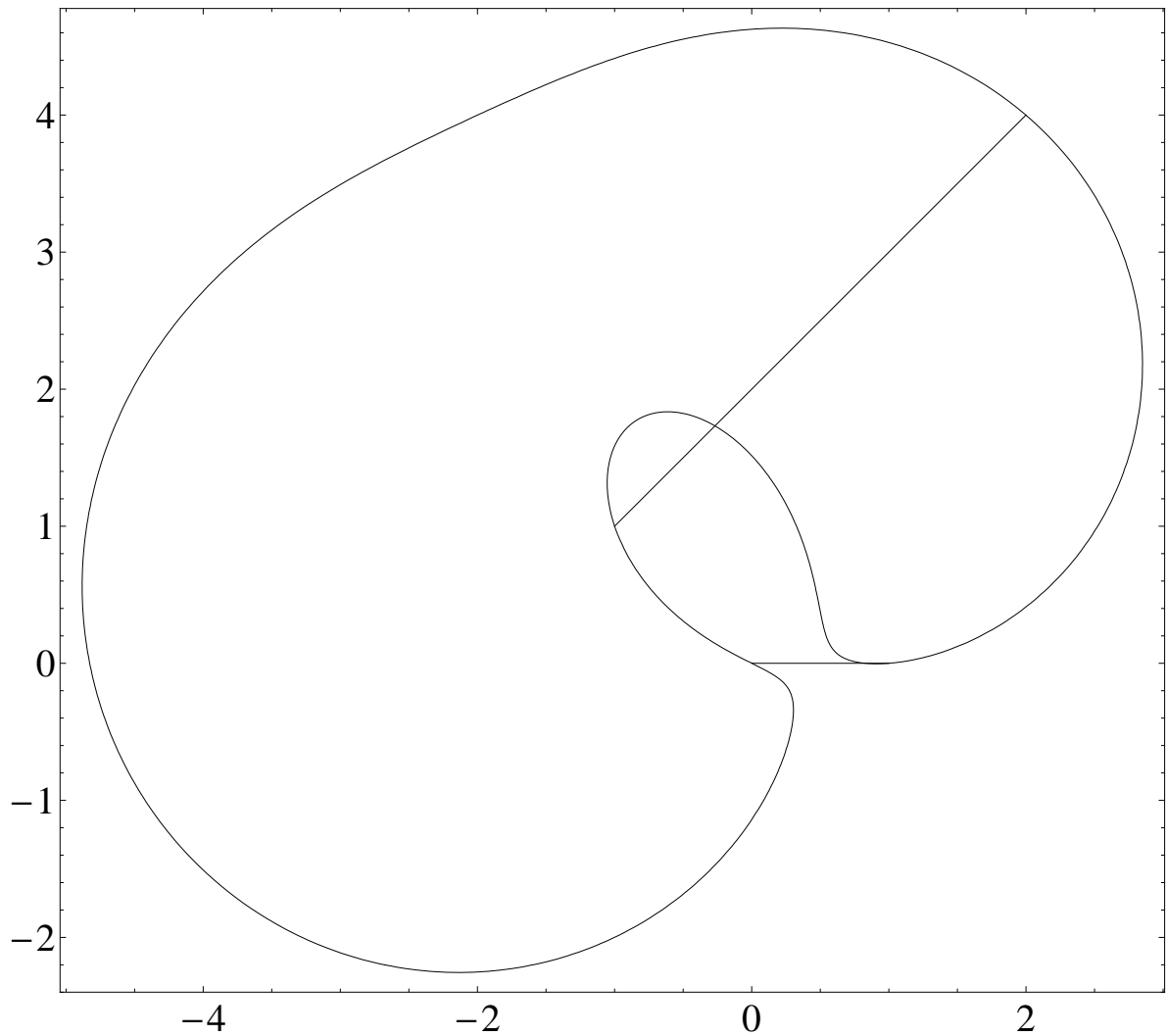}
\caption{\footnotesize{The curve $G_{TUVW}=0$ за $T(0,0)$, $U(1,0)$, $V(2,4)$ and $W(-1,1)$. The segments (i.e.~the fixed one-dimensional simplexes) $TU$ and $VW$ are also given.} }
        \label{GTUVW1}
\end{figure}

\begin{figure}[h]
    \includegraphics[width=0.34\textwidth]{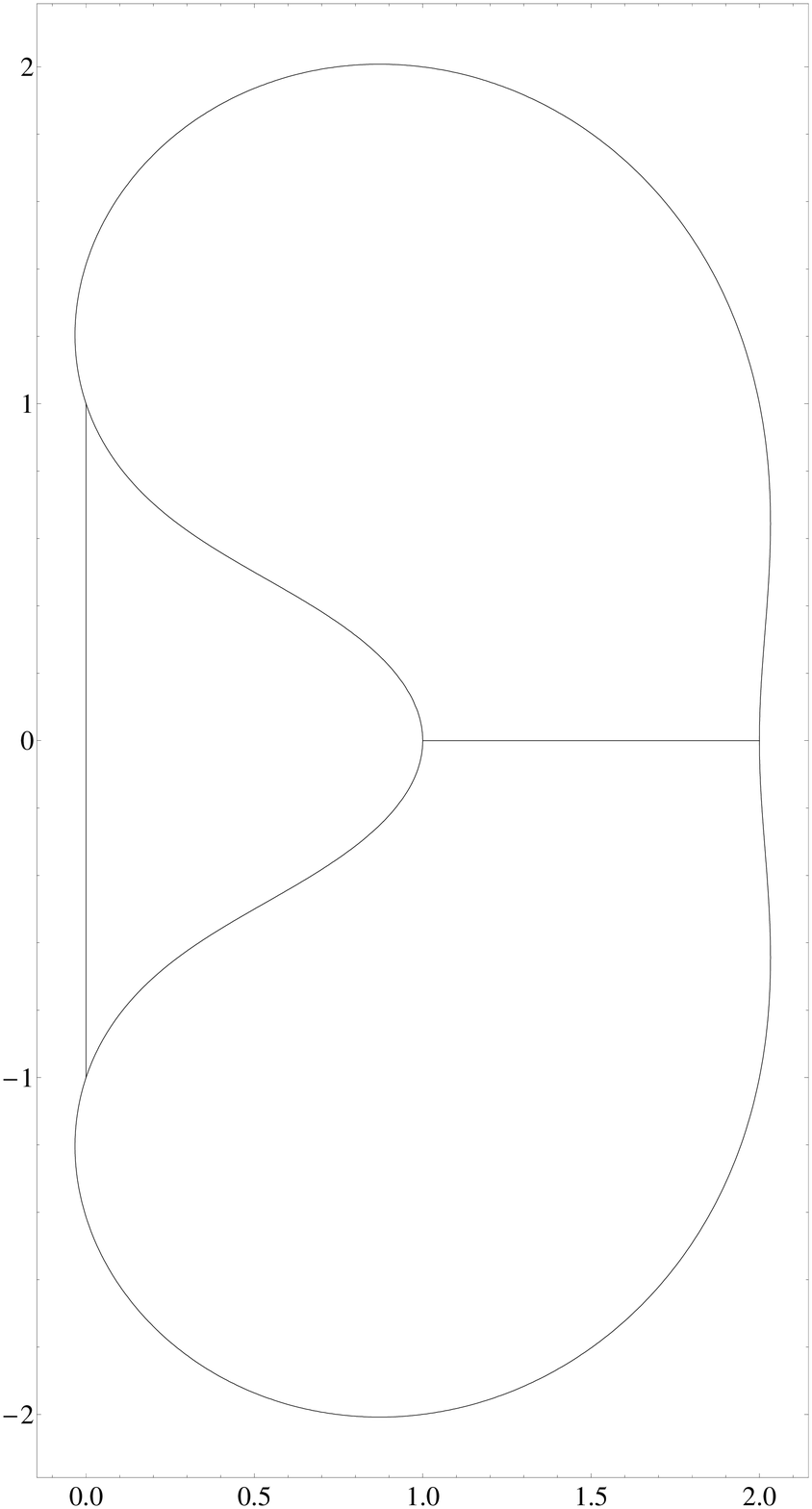}\hfill
    \includegraphics[width=0.64\textwidth]{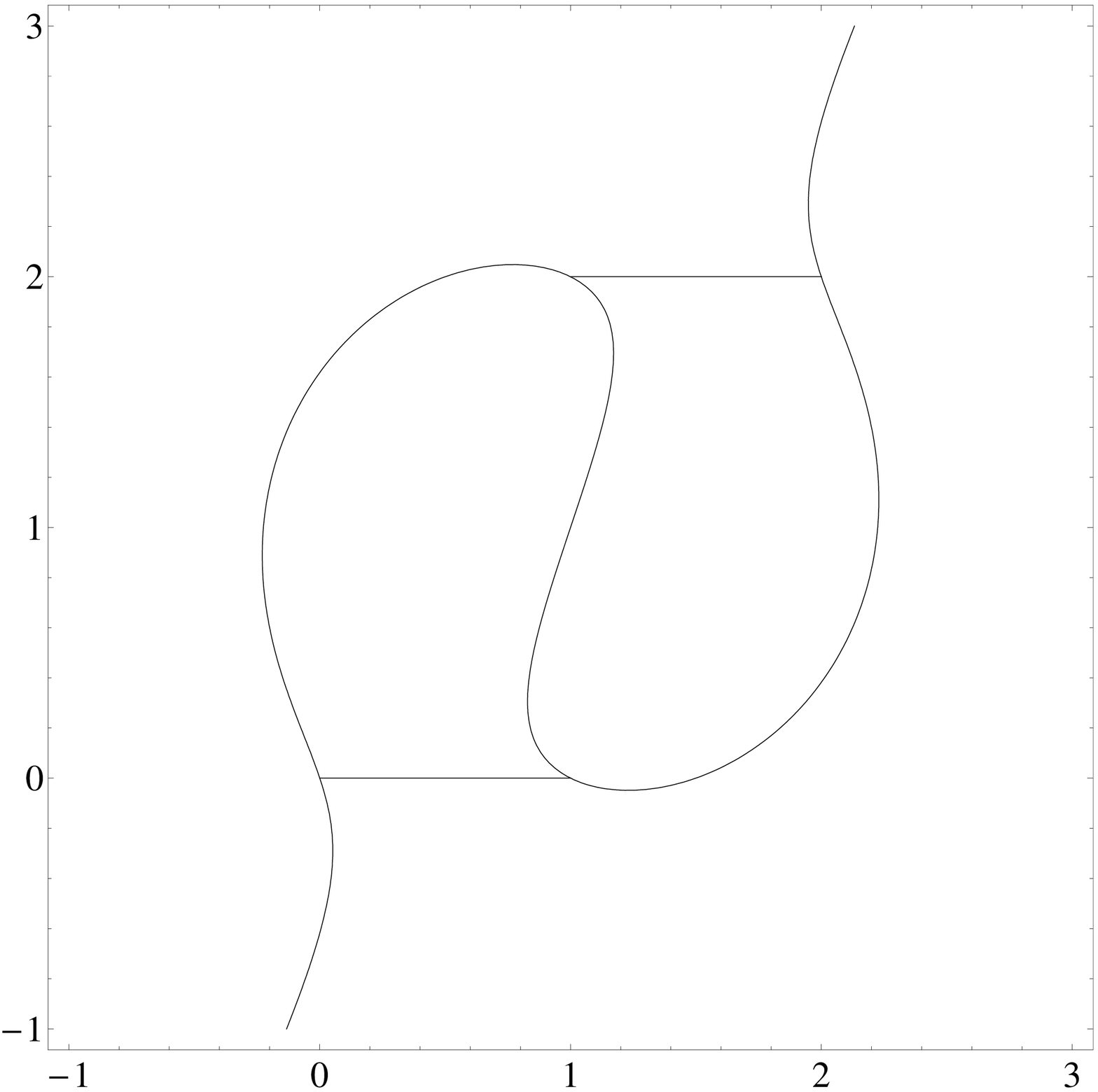}
\caption{\footnotesize{On the left: the curve $G_{TUVW}=0$ for $T(0,1)$, $U(0,-1)$, $V(1,0)$ and $W(2,0)$. On the right: the curve $G_{TUVW}=0$ in the open case. $T(0,0)$, $U(1,0)$, $V(2,2)$ and $W(1,2)$. The segments $TU$ and $VW$ are also given.} }
        \label{GTUVW2}
\end{figure}

\begin{figure}[h]
    \includegraphics[width=0.49\textwidth]{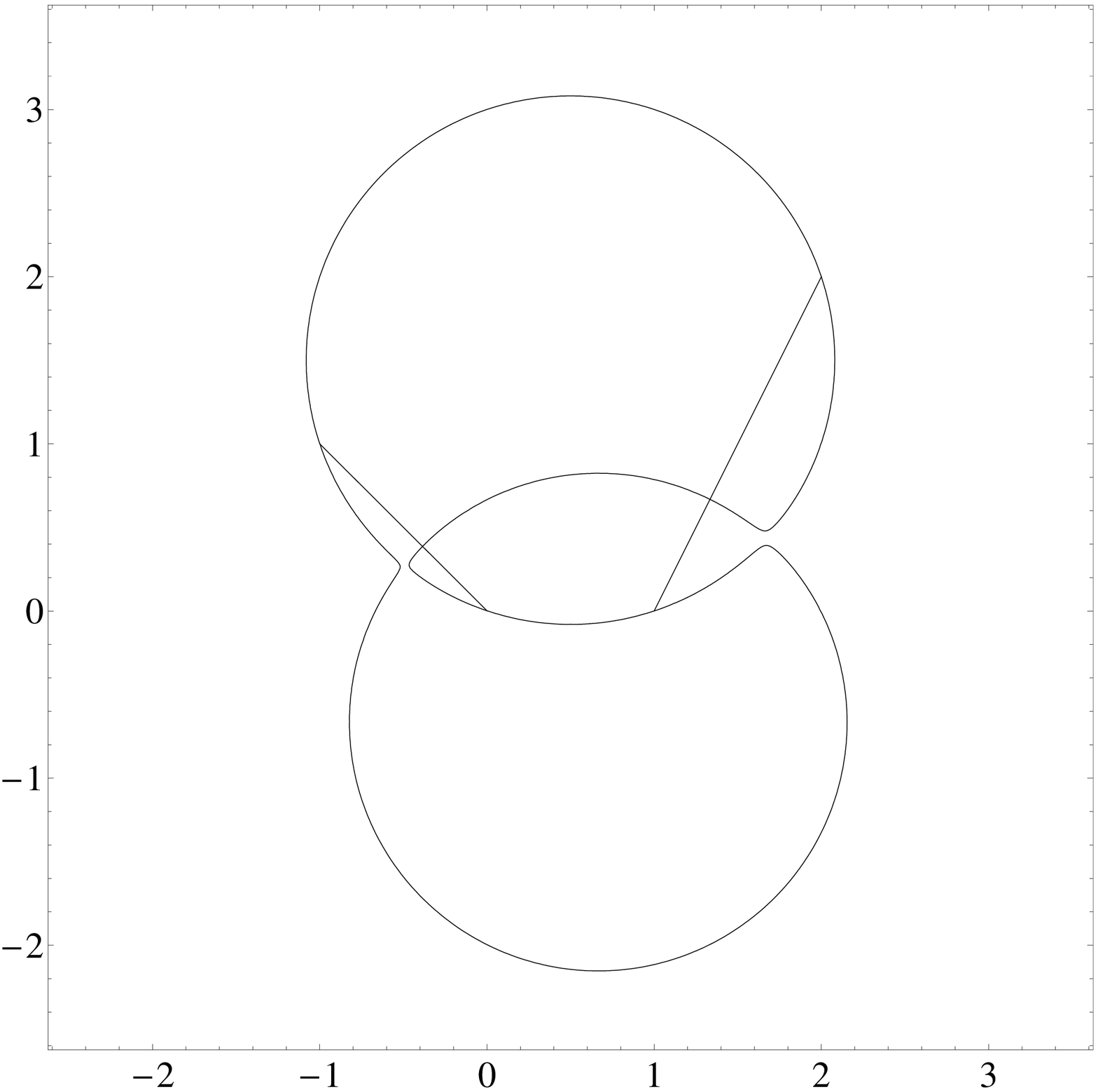}\hfill
    \includegraphics[width=0.49\textwidth]{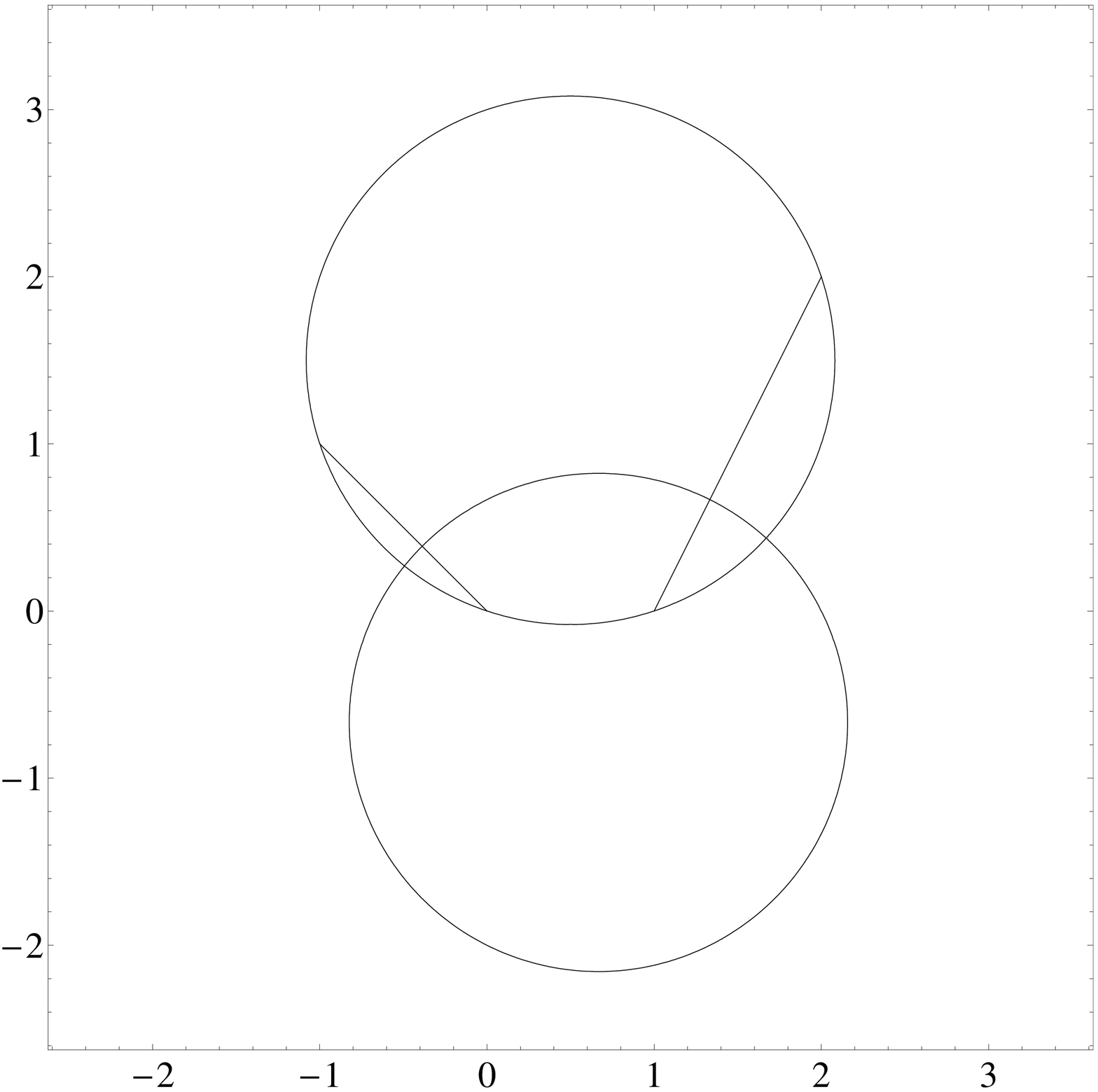}
\caption{\footnotesize{The curve $G_{TUVW}=0$ for $T(1,0)$, $U(2,2)$ and $V(-1,1)$. On the left: the non-degenerate case $W(0.001,0.001)$. On the right: the degenerate case, in which it is possible to draw a circle around $TUVW$---for $W(0,0)$. The segments $TU$ and $VW$ are also given.} }
        \label{GTUVW4}
\end{figure}

\begin{figure}[h]
    \includegraphics[width=0.49\textwidth]{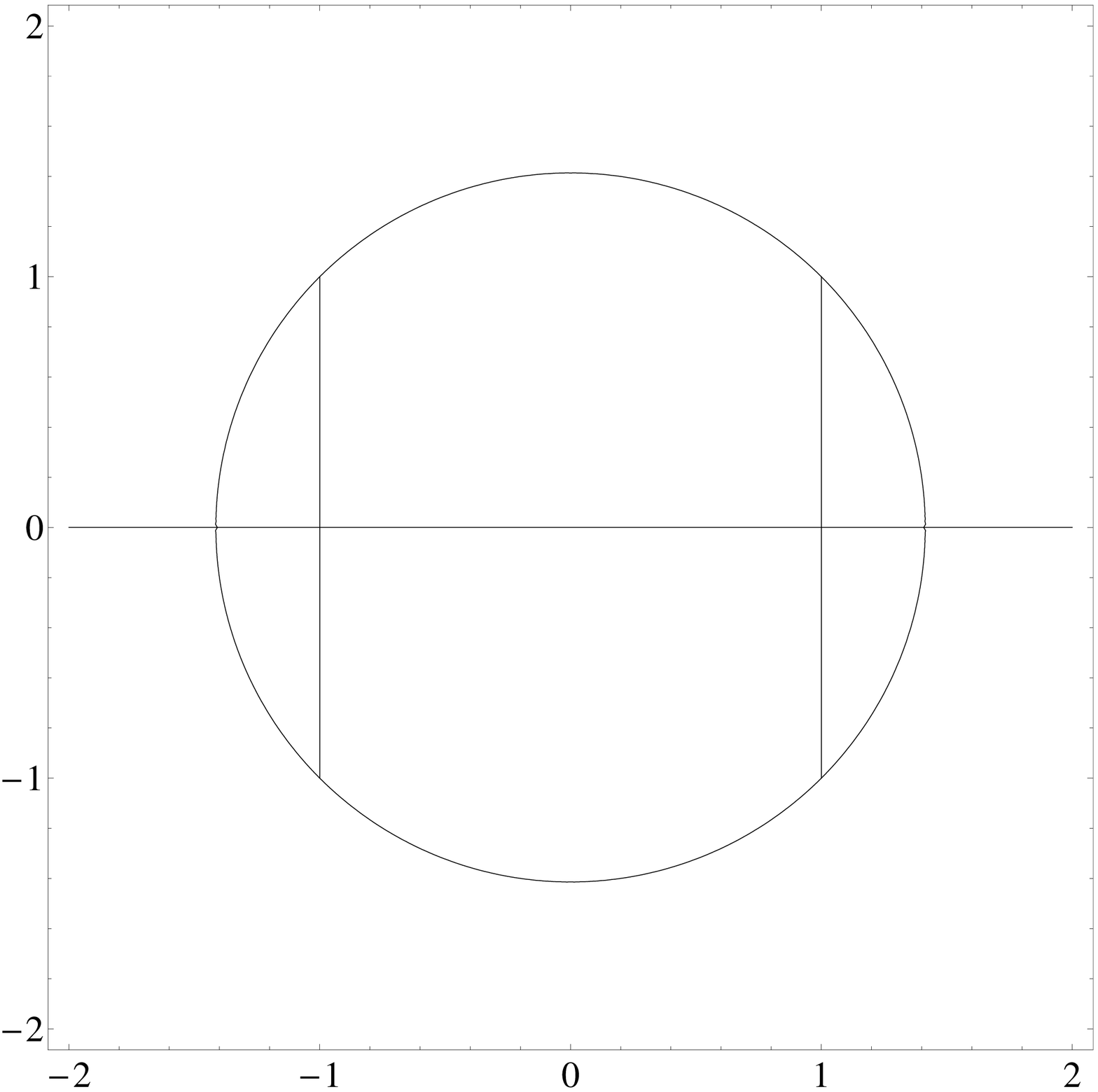}\hfill \includegraphics[width=0.49\textwidth]{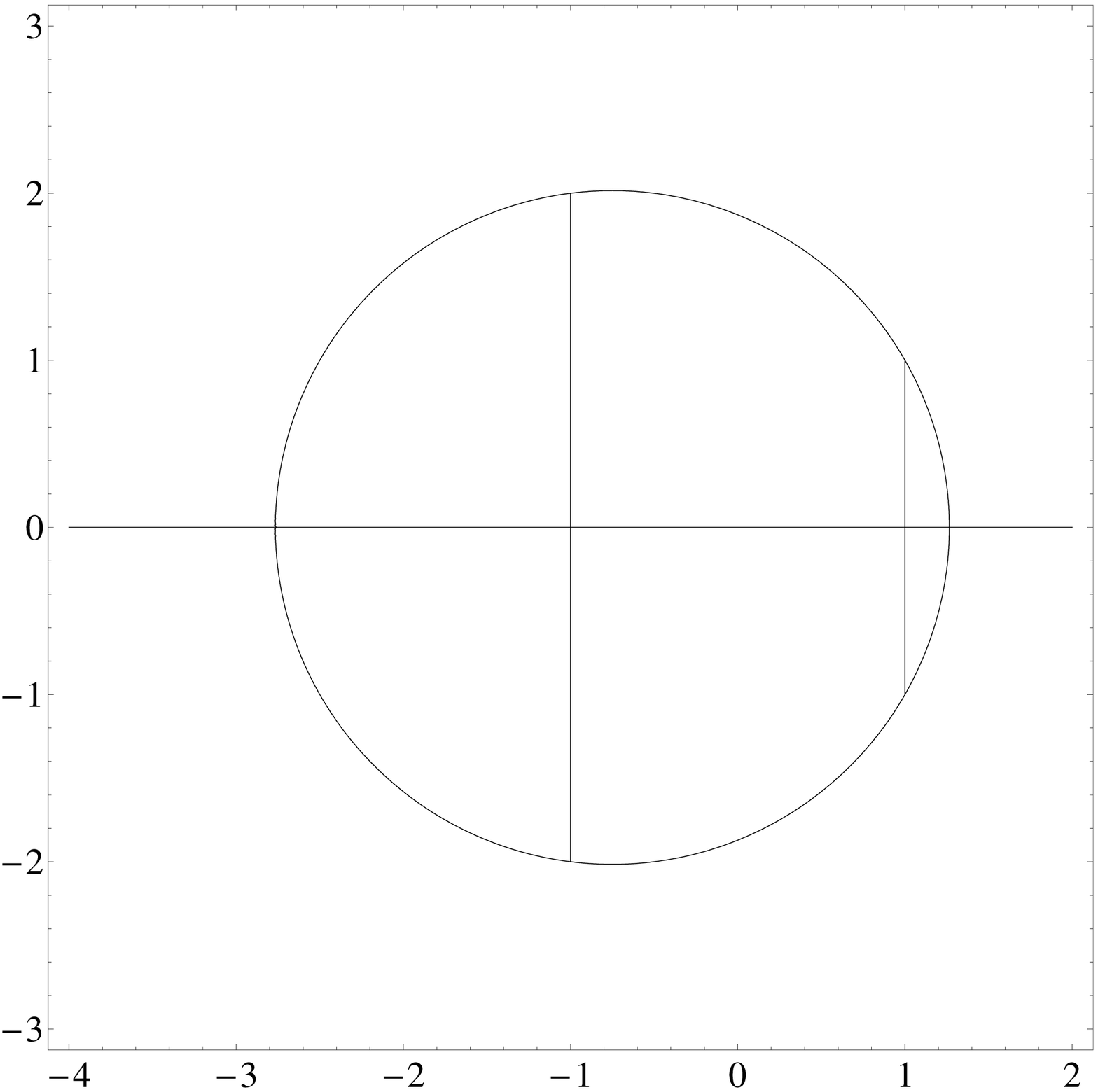}
\caption{\footnotesize{Graphs of the curves $G_{TUVW}=0$ in the degenerate cases. On the left: the case when $TUVW$ is a rectangle. $T(-1,1)$, $U(-1,-1)$, $V(1,-1)$ and $W(1,1)$. On the right: the case when $TUVW$ is a bilateral trapezium  ($TU\parallel VW$). $T(-1,2)$, $U(-1,-2)$, $V(1,-1)$ and $W(1,1)$. The segments $TU$ and $VW$ are also given.} }
        \label{GTUVW6}
\end{figure}

\begin{figure}[h]
    \includegraphics[width=0.49\textwidth]{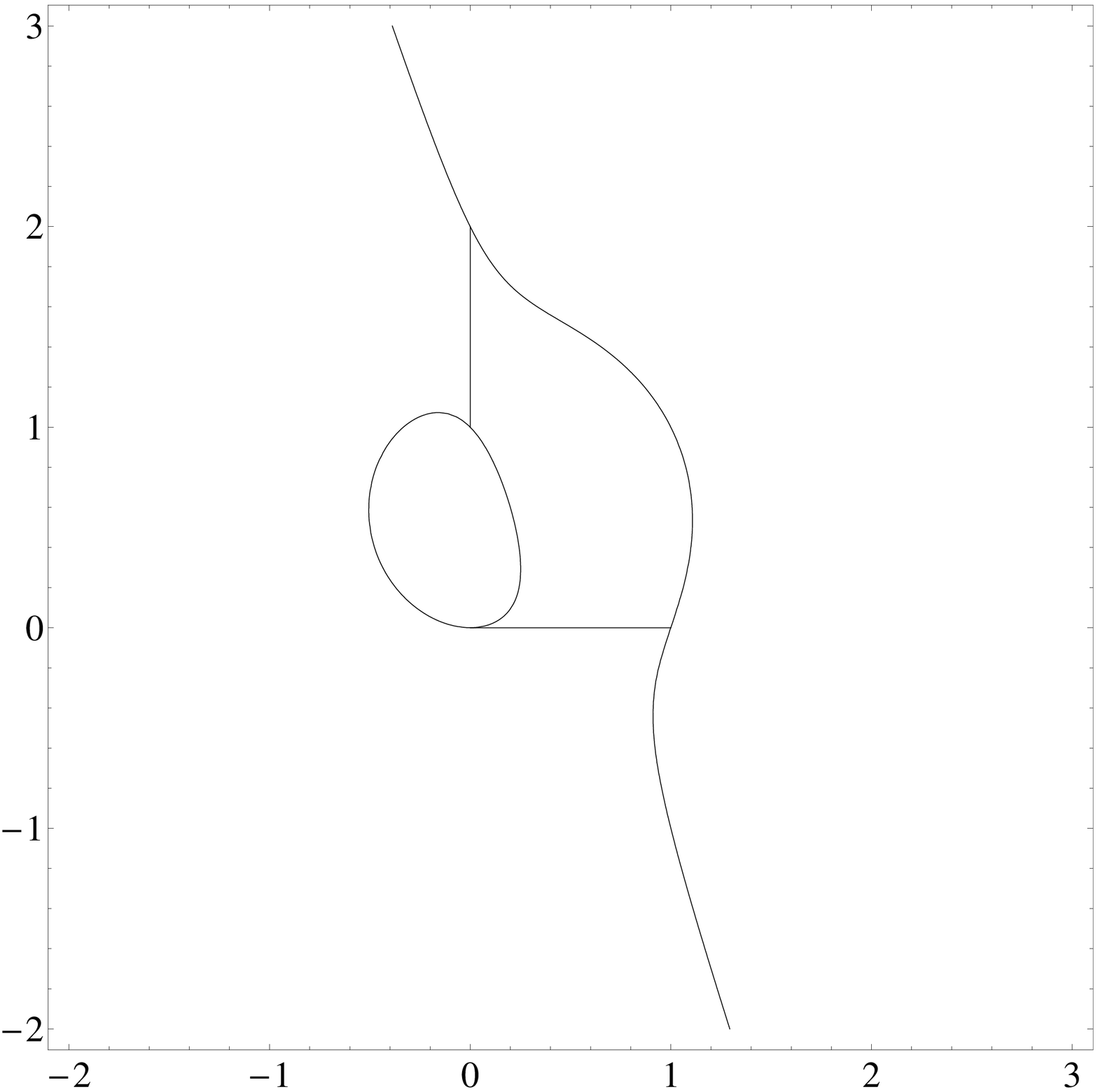}\hfill
    \includegraphics[width=0.49\textwidth]{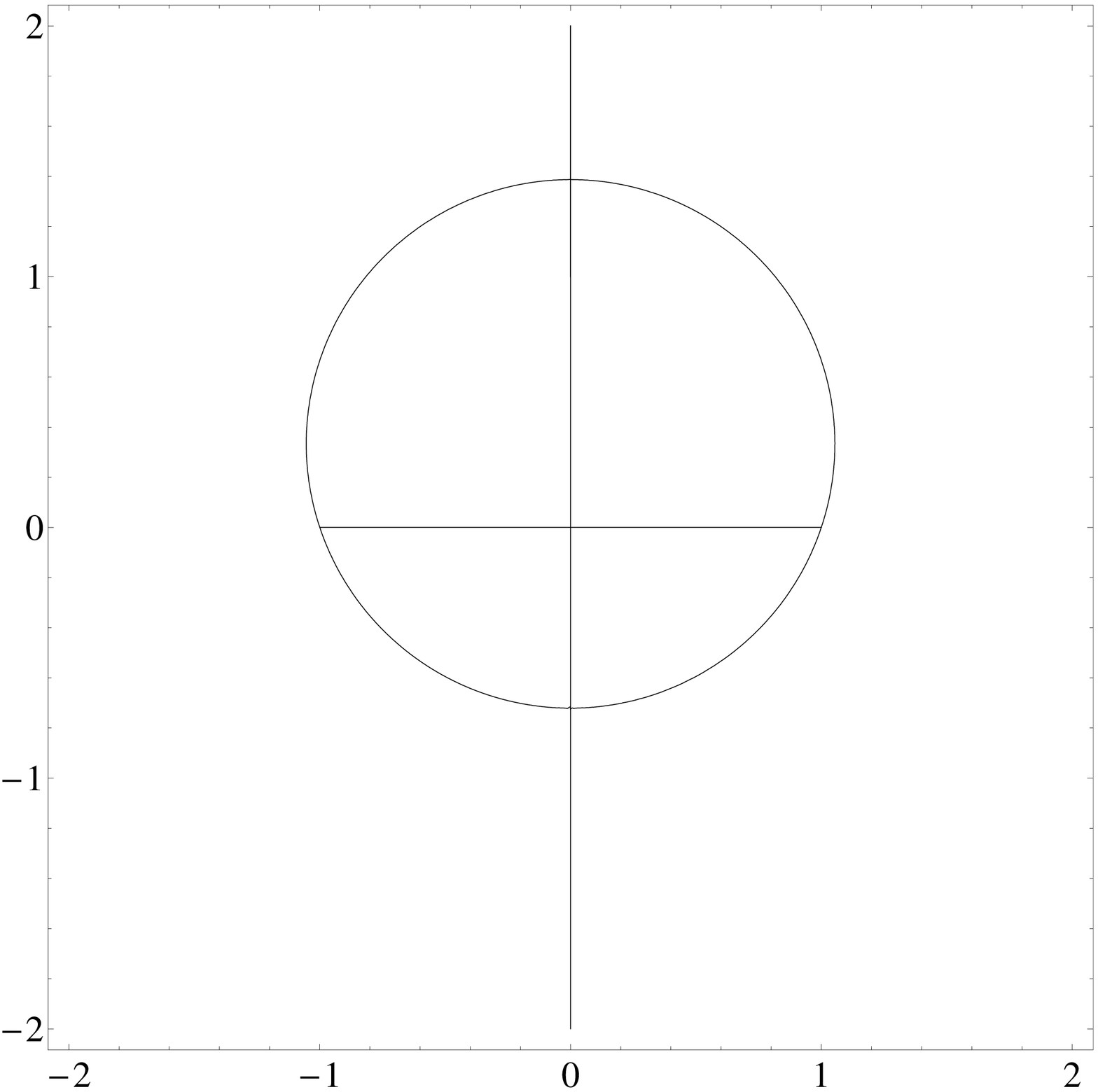}
\caption{\footnotesize{The curve $H_{TUVW}=0$ for $\angle(TU,VW)=\pi/2$ rad. On the left: $T(0,0)$, $U(1,0)$, $V(0,1)$ and $W(0,2)$; on the right: the degenerate case when $VW$ lies on the perpendicular bisector of $TU$, $T(-1,0)$, $U(1,0)$, $V(0,1)$ and $W(0,2)$. The segments $TU$ and $VW$ are also given.} }
        \label{GTUVW7a}
\end{figure}

\begin{figure}[h]
    \includegraphics[width=0.98\textwidth]{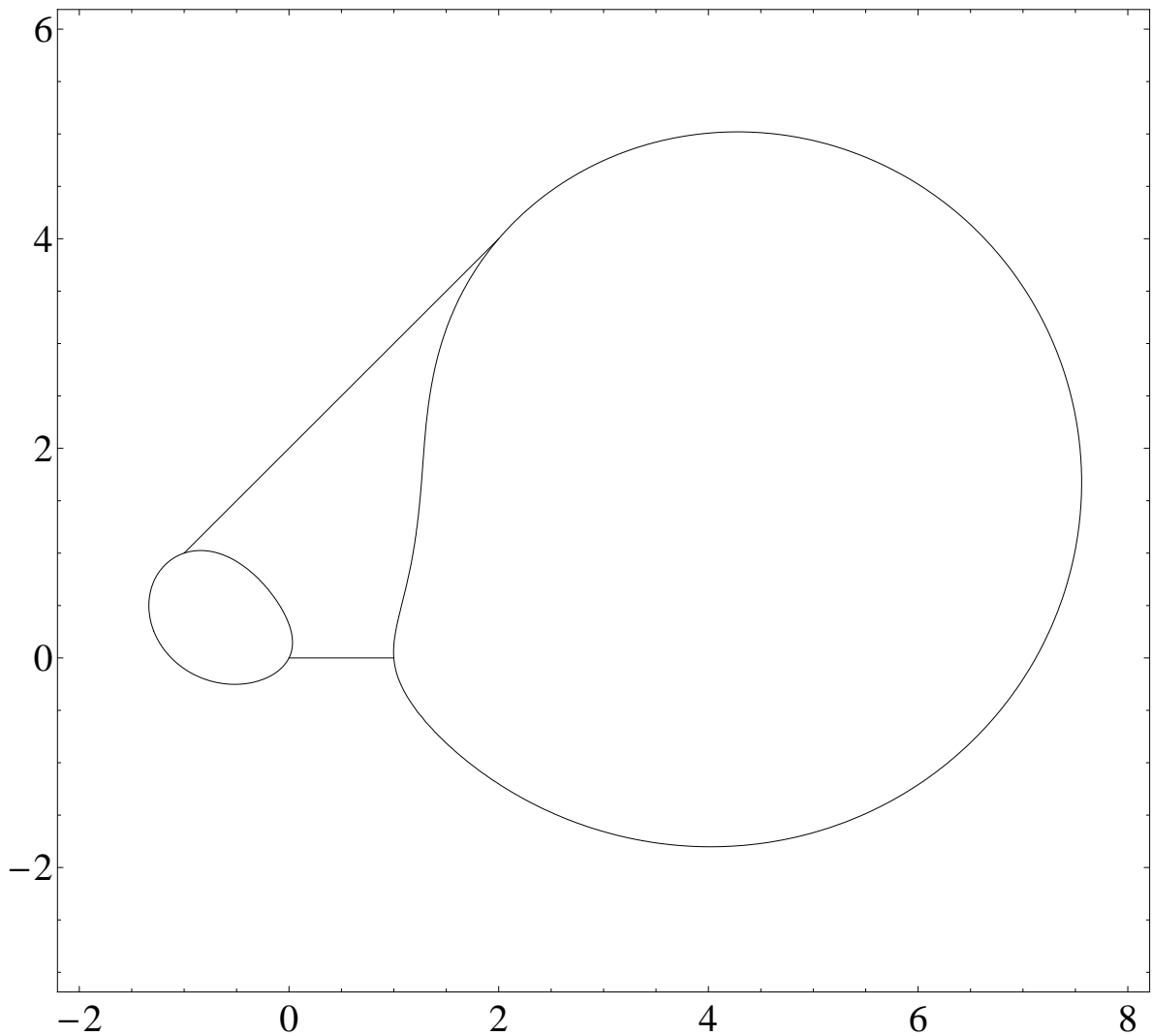}
\caption{\footnotesize{The curve $H_{TUVW}=0$ for $T(0,0)$, $U(1,0)$, $V(2,4)$ and $W(-1,1)$. The segments $TU$ and $VW$ are also given.} }
        \label{GTUVW7}
\end{figure}

\begin{figure}[h]
    \includegraphics[width=0.49\textwidth]{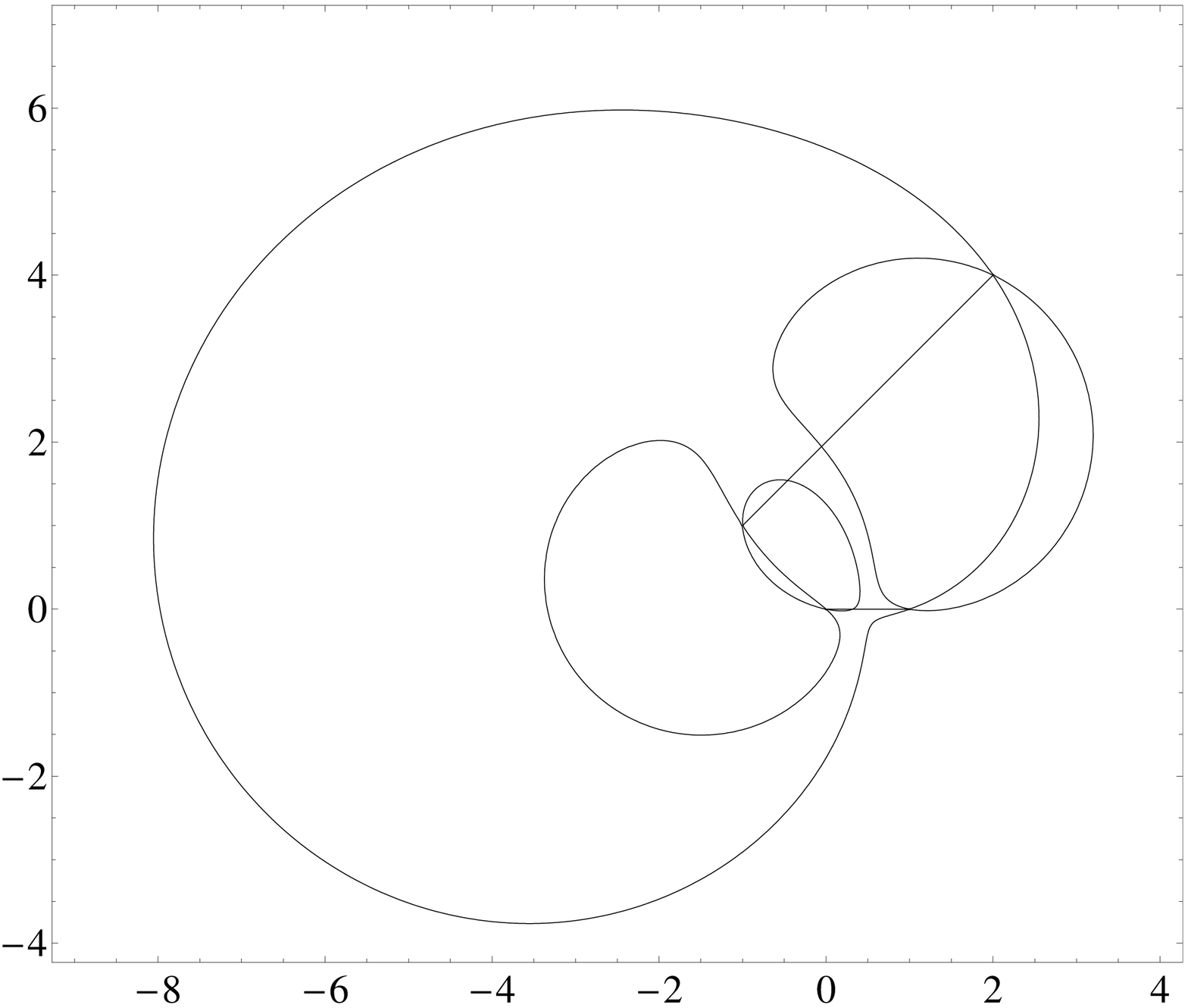}\hfill
    \includegraphics[width=0.49\textwidth]{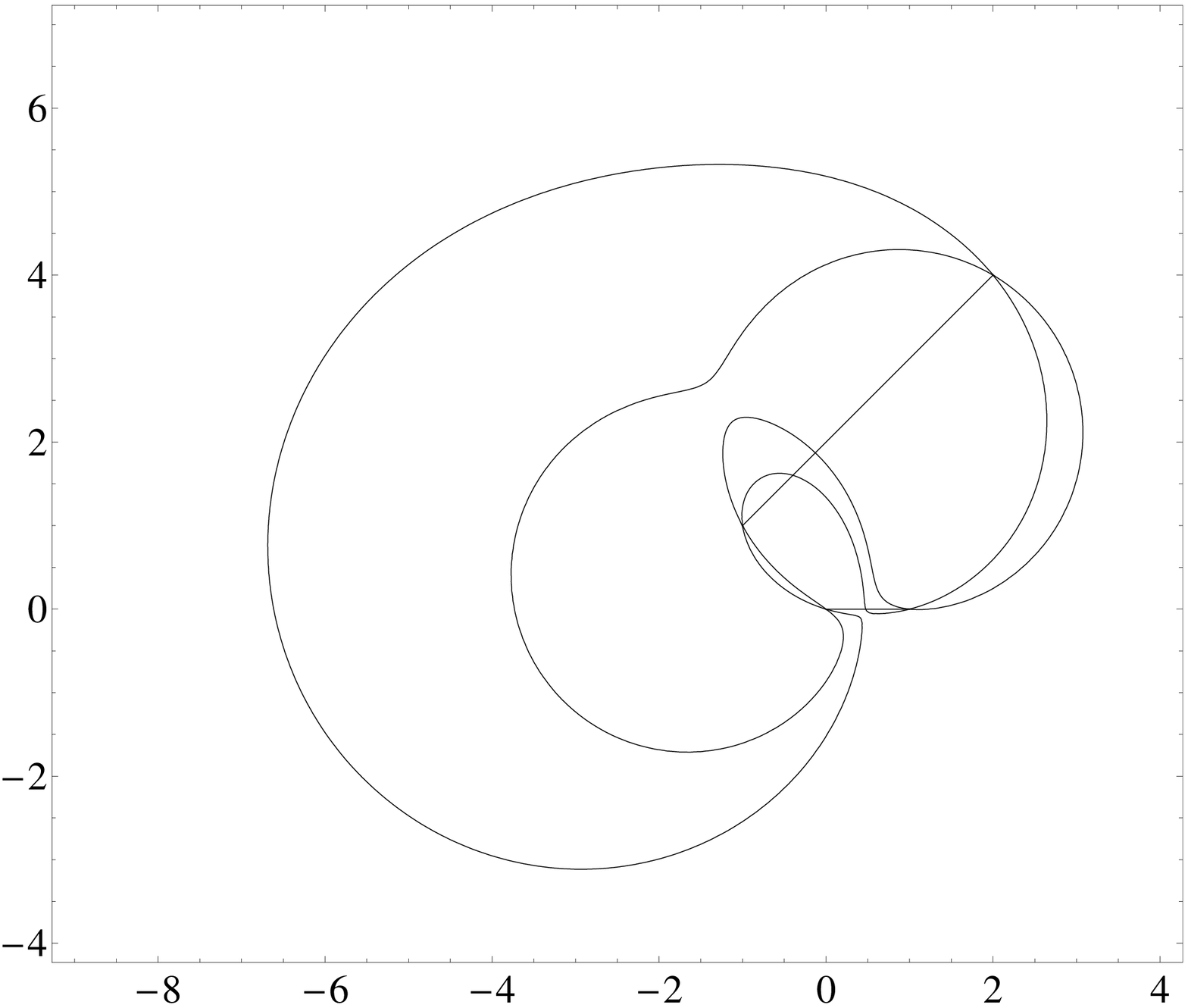}
\caption{\footnotesize{The curve $F_{TUVW}=0$ for $T(0,0)$, $U(1,0)$, $V(2,4)$ and $W(-1,1)$. On the left: $\alpha=\pi/12$ rad; on the right: $\alpha=\pi/18$ rad. The segments $TU$ and $VW$ are also given.} }
        \label{GTUVW8}
\end{figure}

\begin{figure}[h]
    \includegraphics[width=0.49\textwidth]{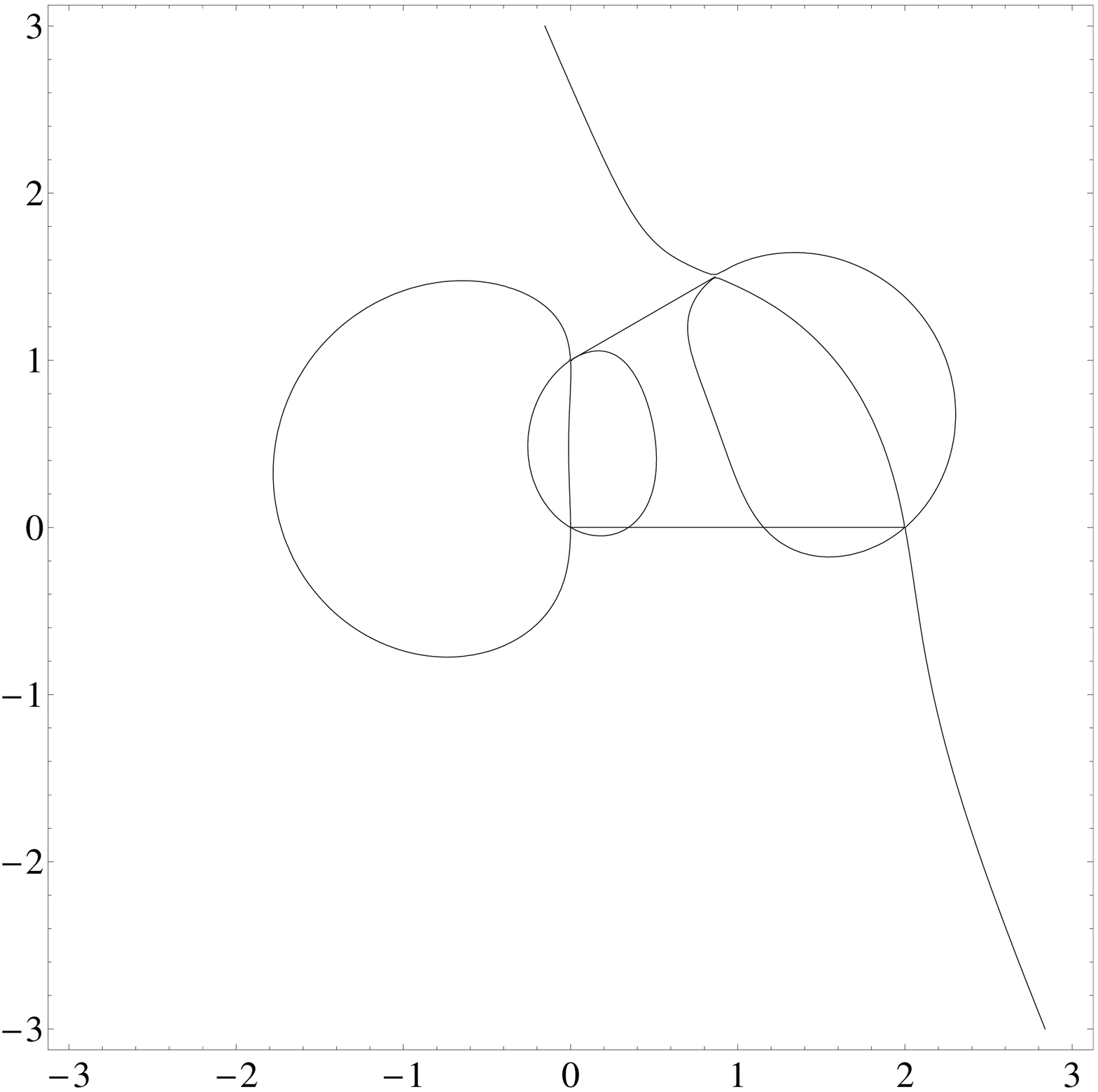}\hfill
    \includegraphics[width=0.49\textwidth]{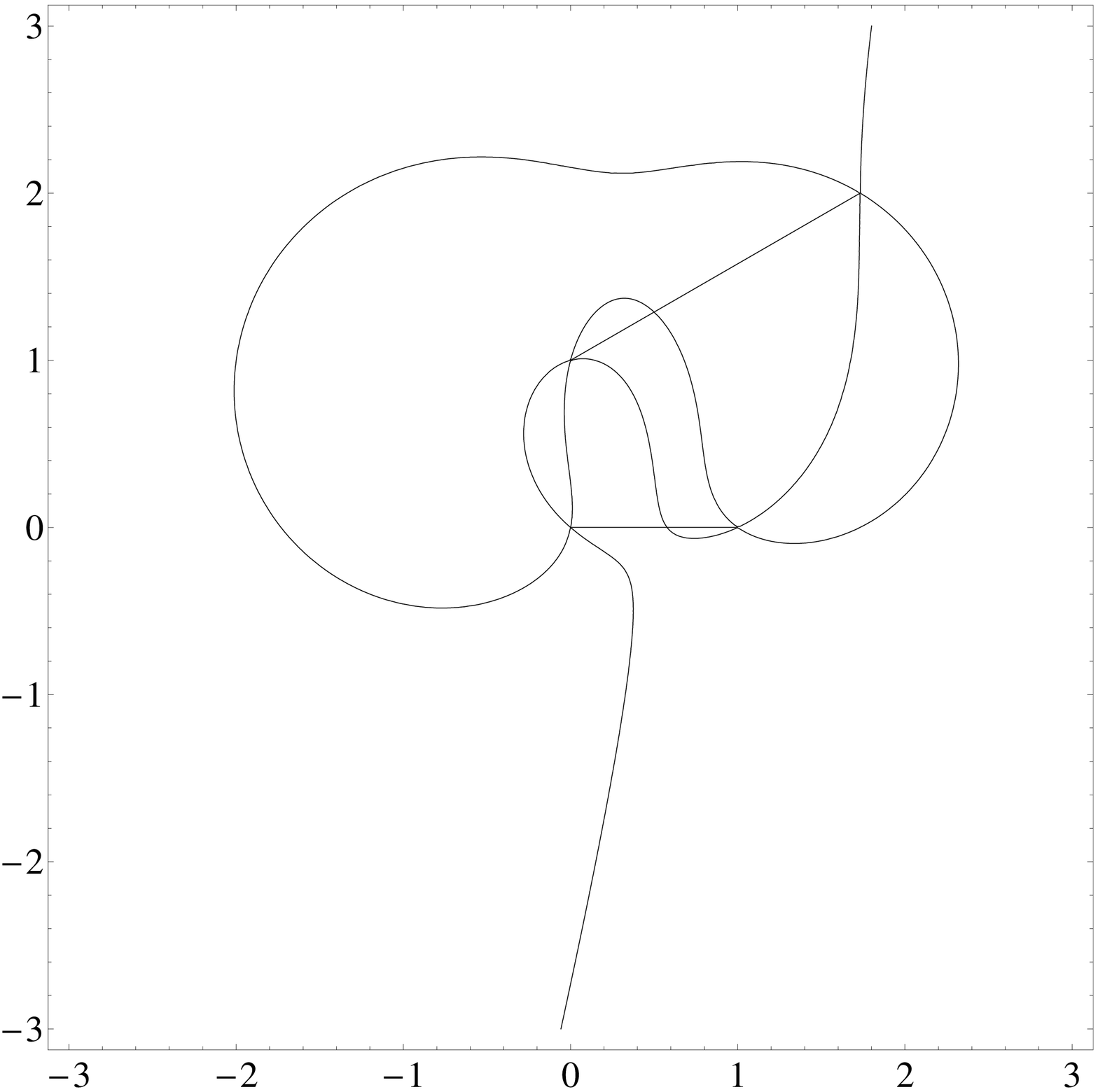}
\caption{\footnotesize{The curve $F_{TUVW}=0$ for $\alpha=\angle(TU,VW)=\pi/6$ rad. On the left: $T(0,0)$, $U(2,0)$, $V(\sqrt{3}/2,3/2)$ and $W(0,1)$; on the right: $T(0,0)$, $U(1,0)$, $V(\sqrt{3},2)$ and $W(0,1)$. The segments $TU$ and $VW$ are also given.} }
        \label{GTUVW8a}
\end{figure}

\begin{figure}[h]
    \includegraphics[width=0.49\textwidth]{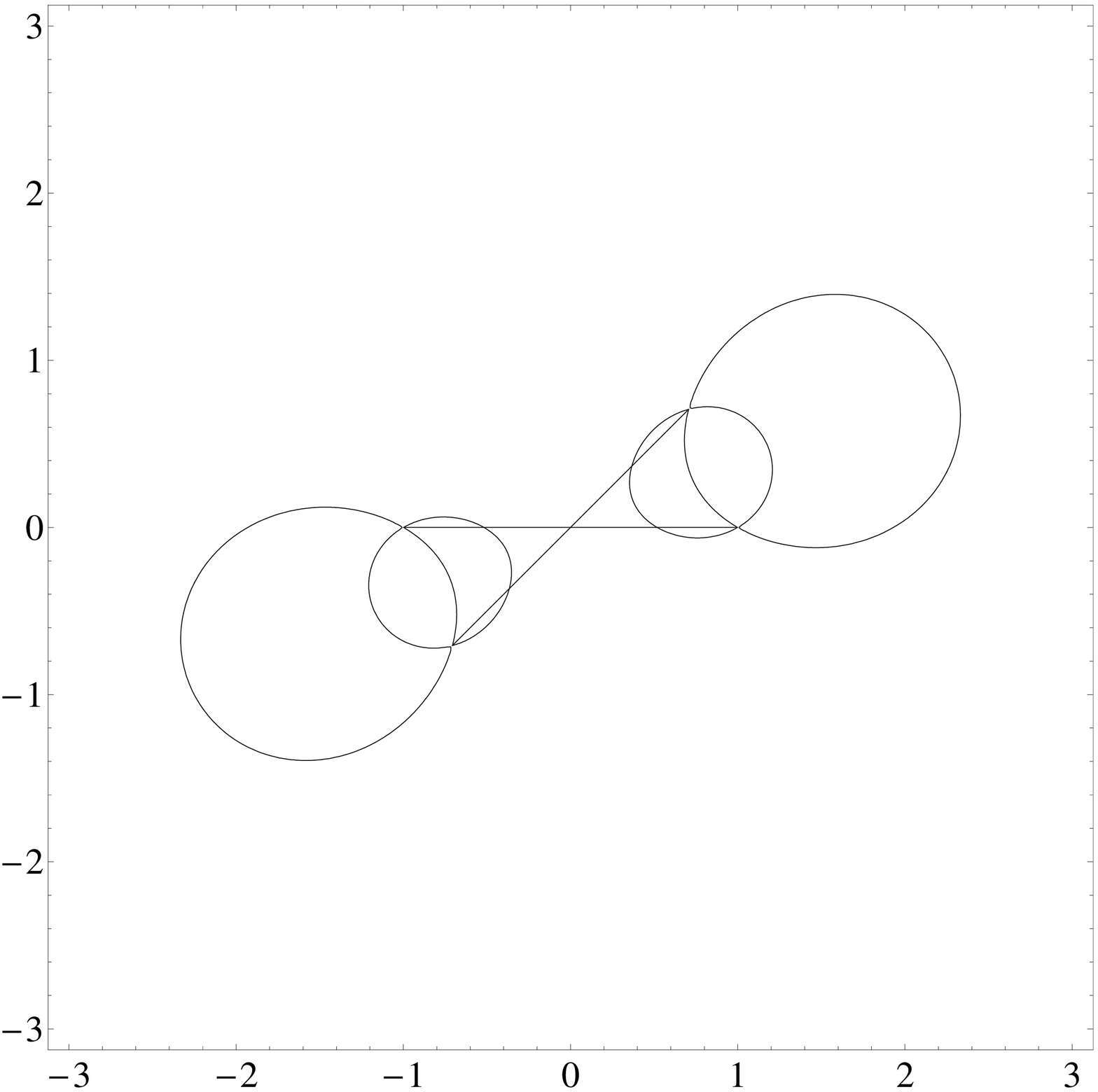}\hfill
    \includegraphics[width=0.49\textwidth]{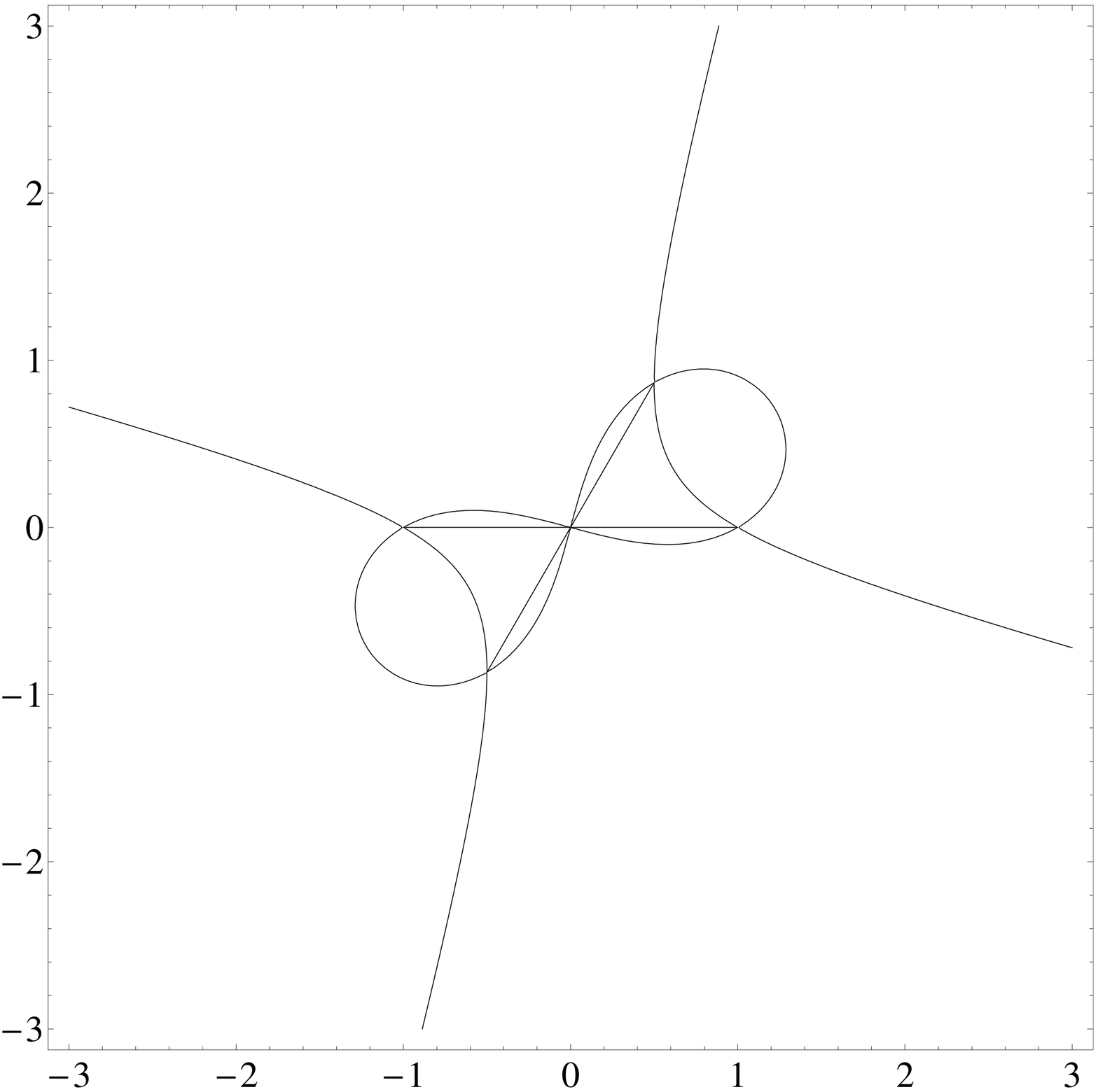}
\caption{\footnotesize{The curve $F_{TUVW}=0$ for $T(-1,0)$, $U(1,0)$, $V(\cos\beta_{12},\sin\beta_{12})$ and $W(-\cos\beta_{12},-\sin\beta_{12})$, where $\beta_{12}=\angle(TU,VW)$. On the left illustration we have $\alpha=\pi/3$ rad and $\beta_{1}=\pi/4$ rad; on the right: $\alpha=\beta_{2}=\pi/3$ rad. The segments $TU$ and $VW$ are also given.} }
        \label{GTUVW8b}
\end{figure}

\begin{figure}[h]
    \includegraphics[width=0.98\textwidth]{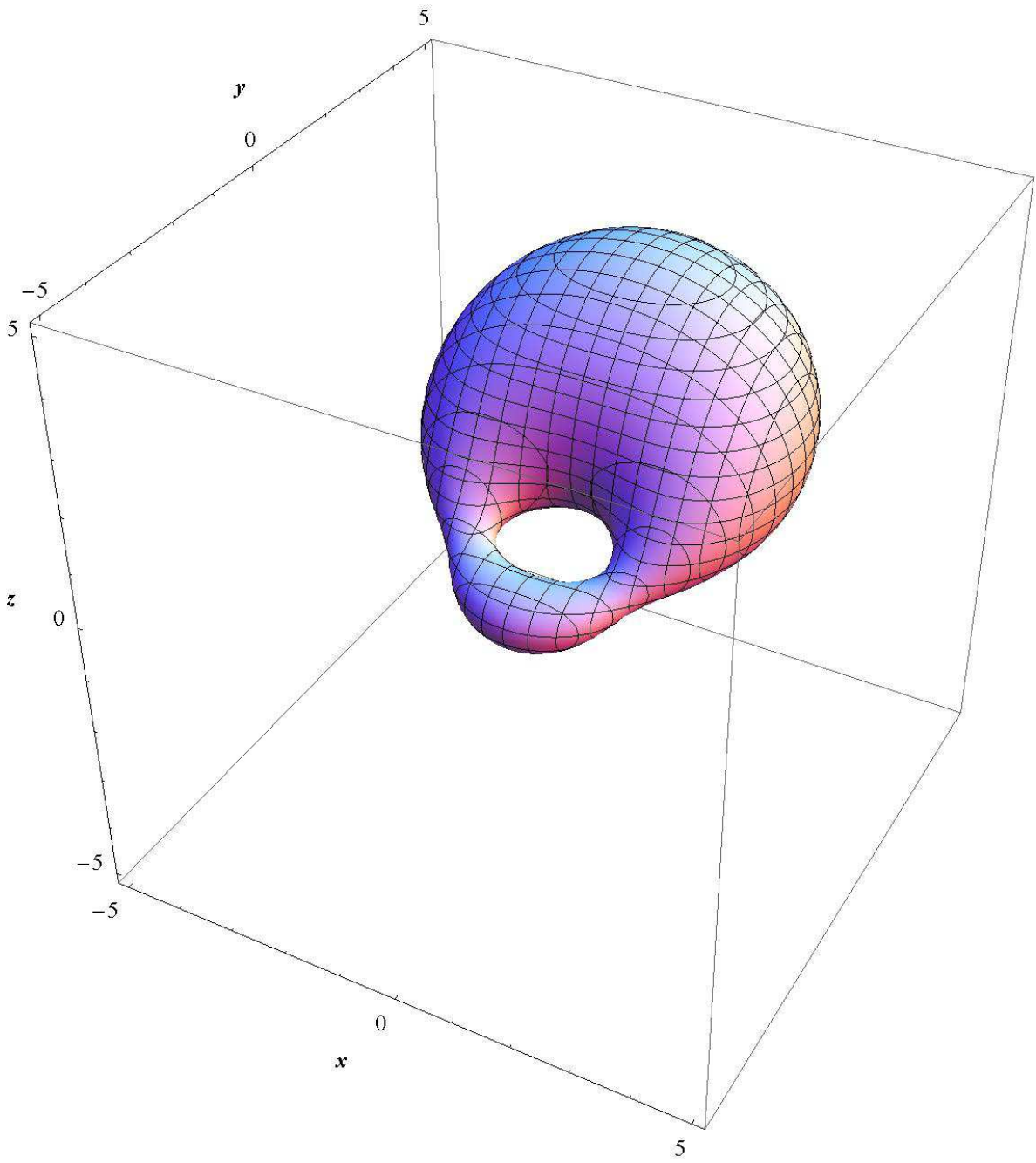}
\caption{\footnotesize{The surface $H_{RSTUVW}=0$ for $R(-1,0,0)$, $S(0,1,0)$, $T(1,0,0)$, $U(-2,0,1)$, $V(0,-1,-1)$ and $W(2,0,1)$.} }
        \label{GTUVW9}
\end{figure}

\begin{figure}[h]
    \includegraphics[width=0.98\textwidth]{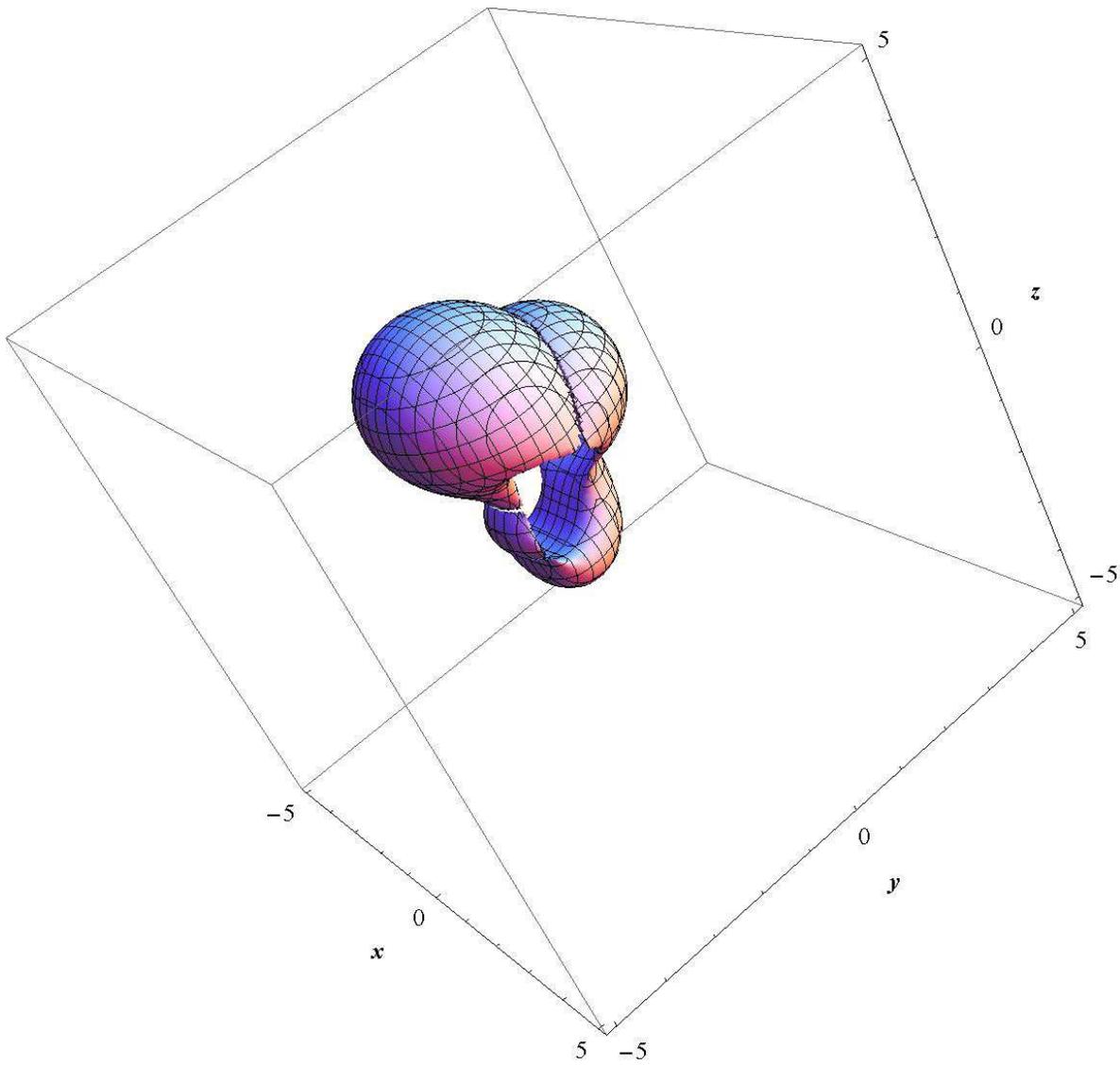}
\caption{\footnotesize{The surface $F_{RSTUVW}=0$ for $R(-1,0,0)$, $S(0,1,0)$, $T(1,0,0)$, $U(-2,0,1)$, $V(0,-1,-1)$ and $W(2,0,1)$ for $\alpha=\pi/6$ rad.} }
        \label{GTUVW10}
\end{figure}

\begin{figure}[h]
        \setlength{\tabcolsep}{ 0 pt }{\scriptsize\tt
        \begin{tabular}{ ccc }
            \includegraphics[width=0.33\textwidth]{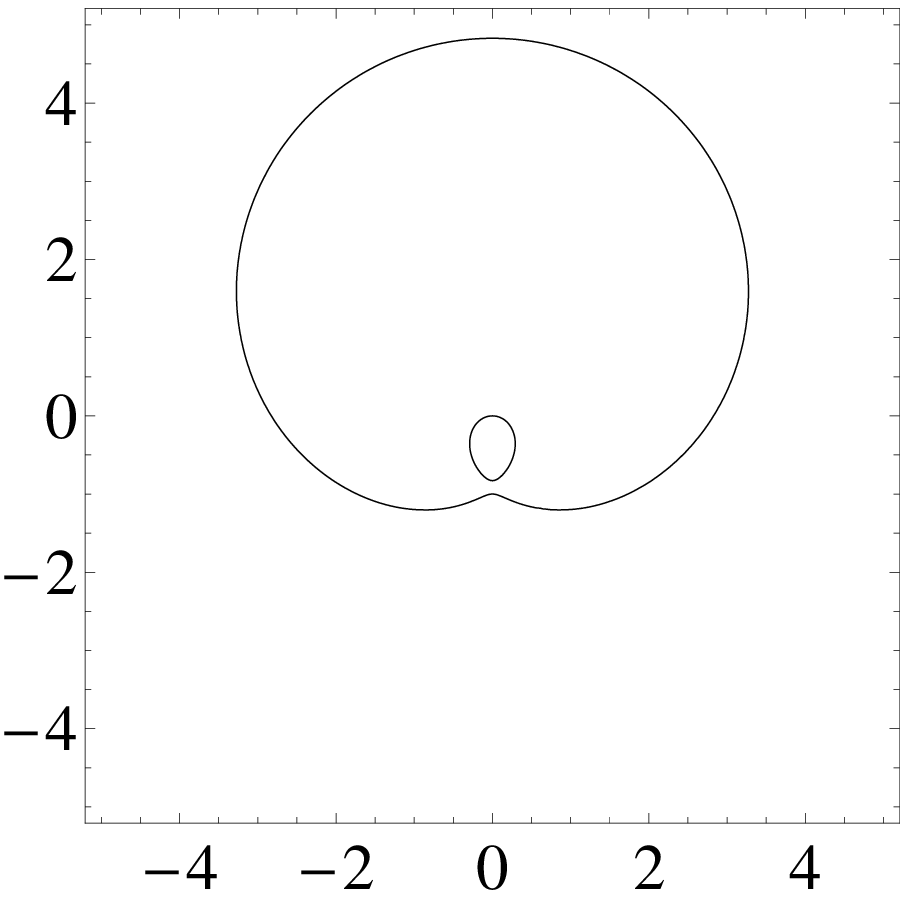} &
            \includegraphics[width=0.33\textwidth]{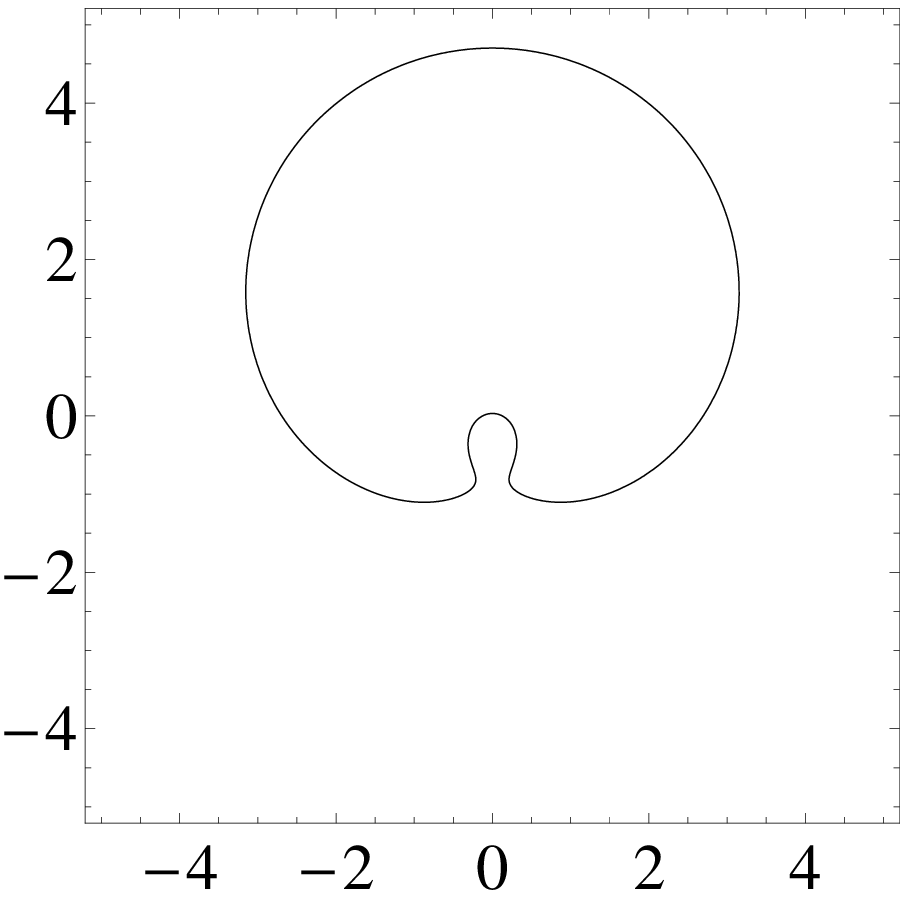} &
            \includegraphics[width=0.33\textwidth]{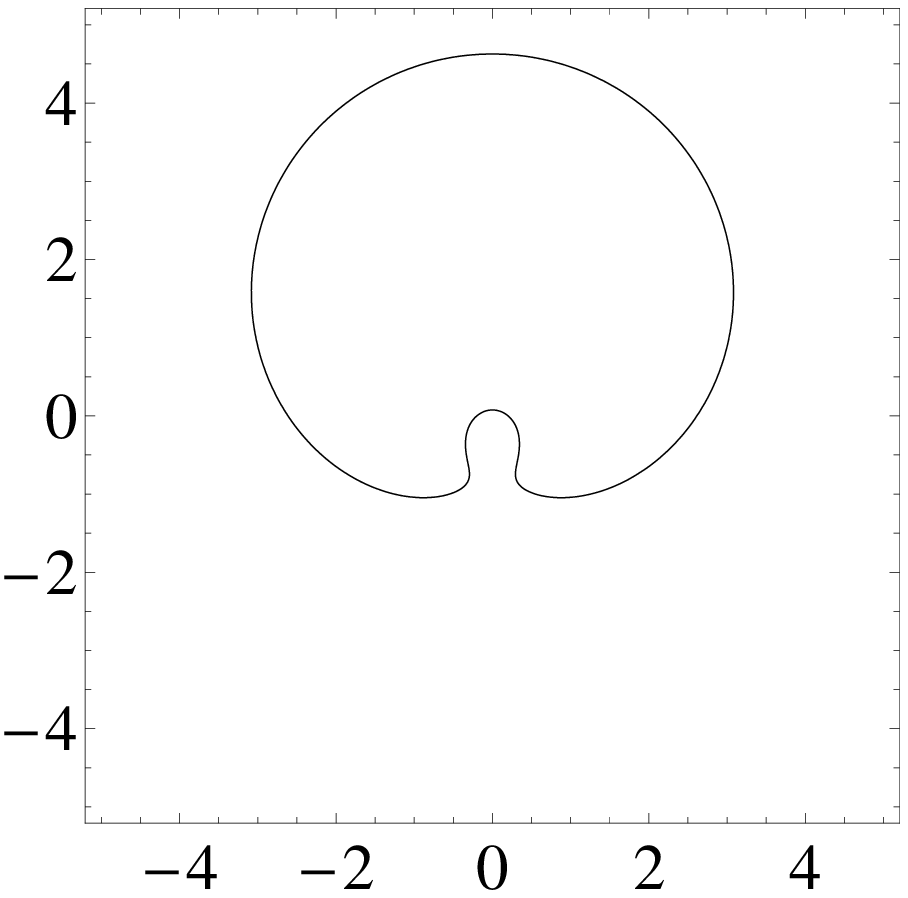} \\
            $t=-1$ &
            $t=-7/8$ &
            $t=-4/5$ \\
            \includegraphics[width=0.33\textwidth]{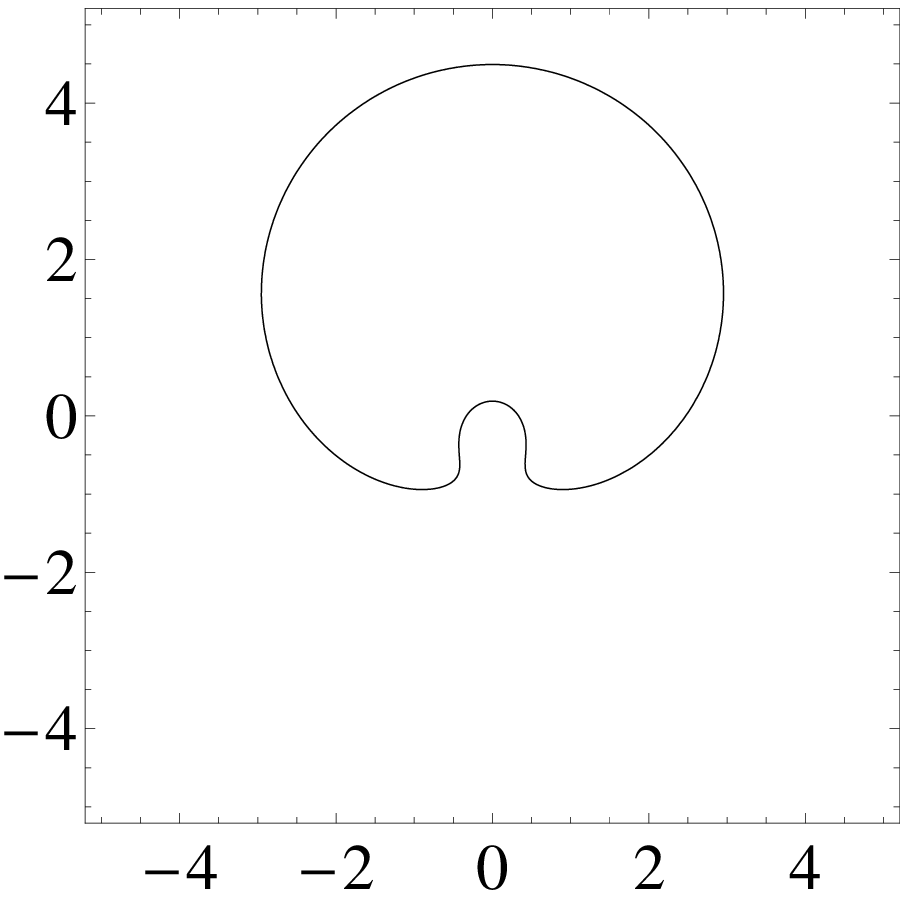} &
            \includegraphics[width=0.33\textwidth]{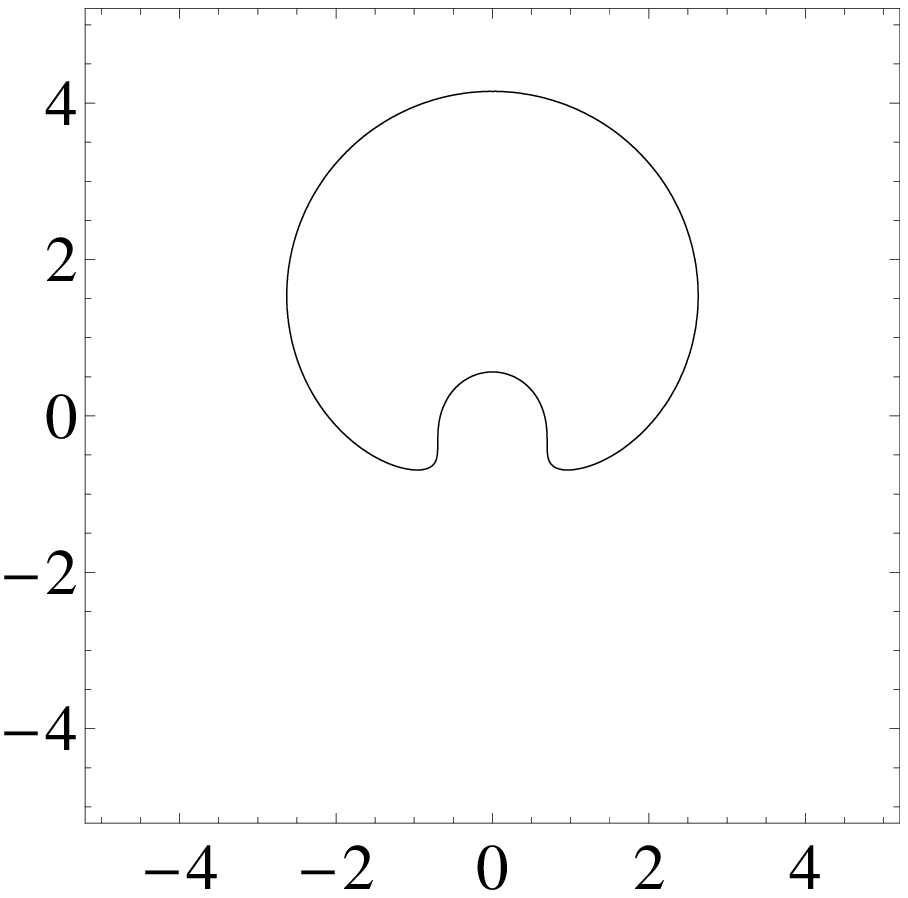} &
            \includegraphics[width=0.33\textwidth]{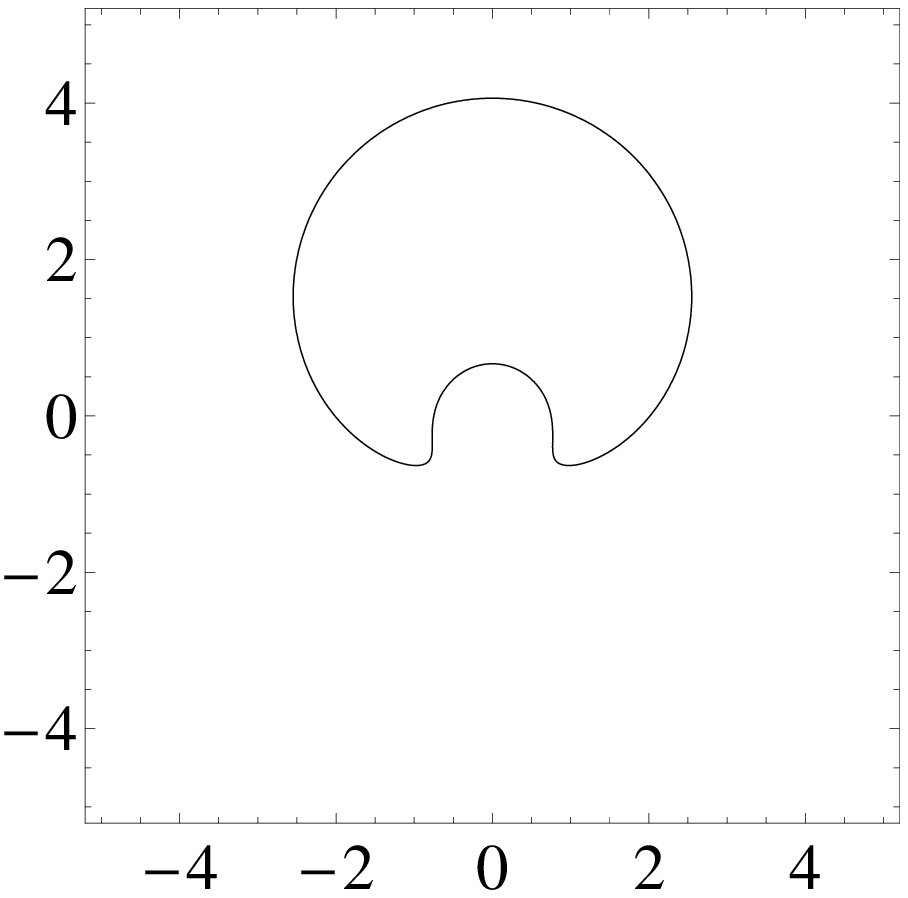} \\
            $t=-2/3$;\  &
            $t=-1/3$;\  &
            $t=-1/4$ \\
            \includegraphics[width=0.33\textwidth]{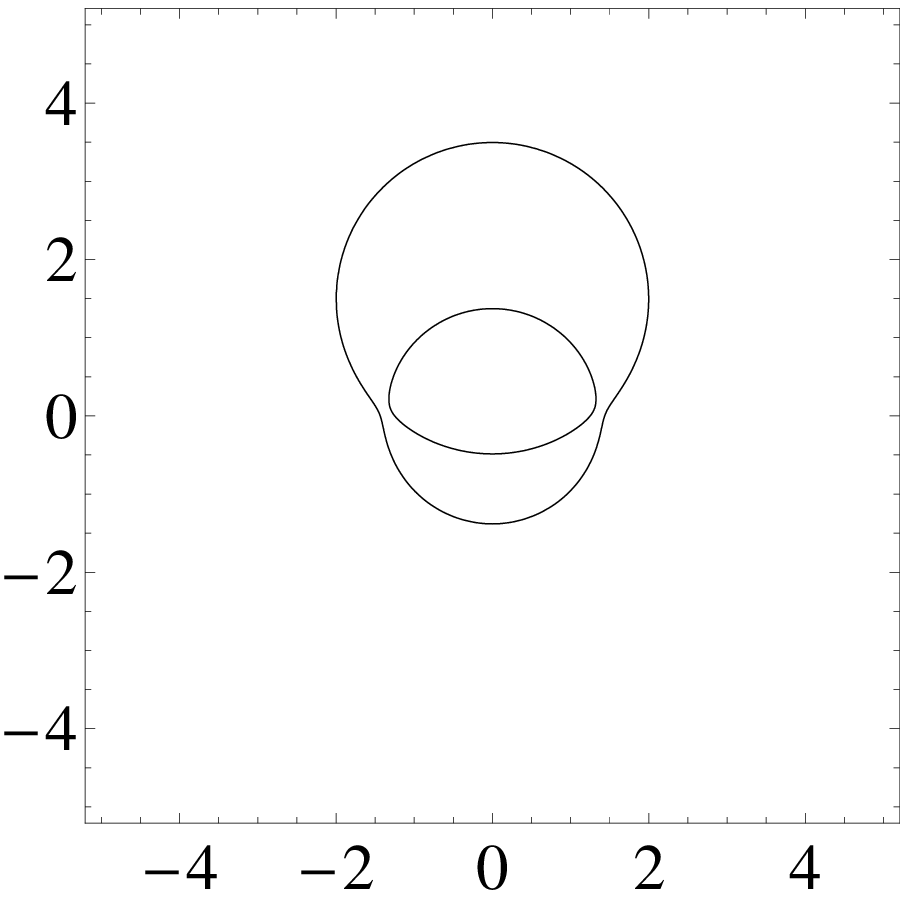} &
            \includegraphics[width=0.33\textwidth]{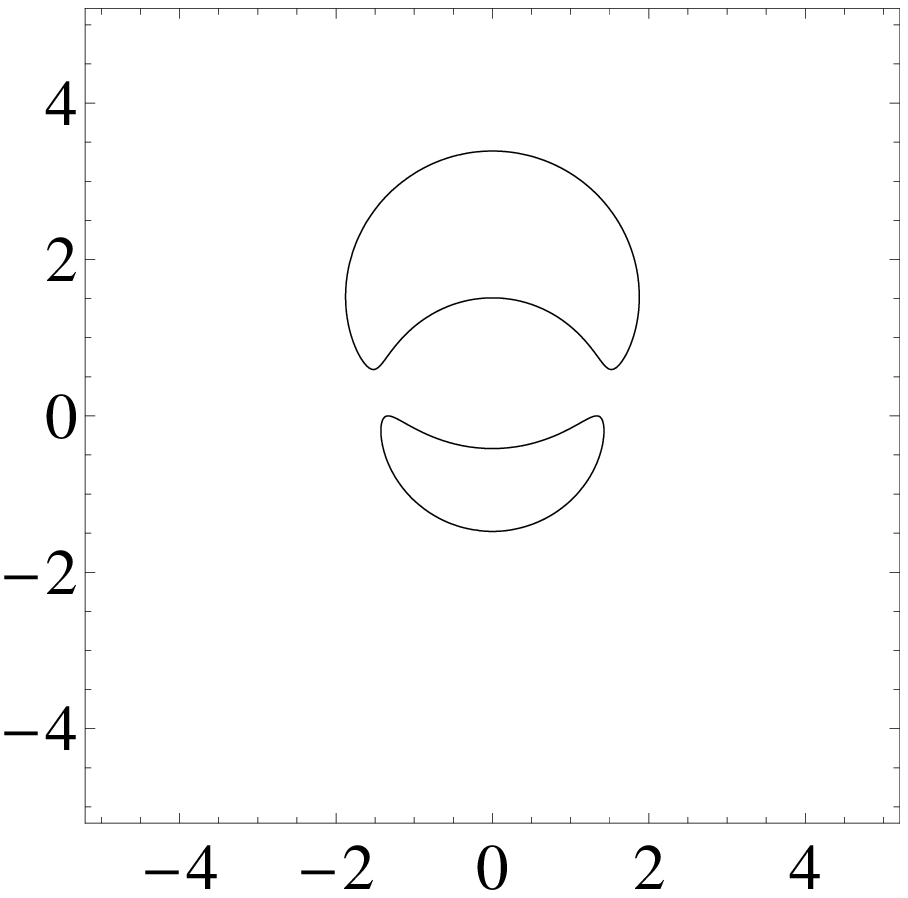} &
            \includegraphics[width=0.33\textwidth]{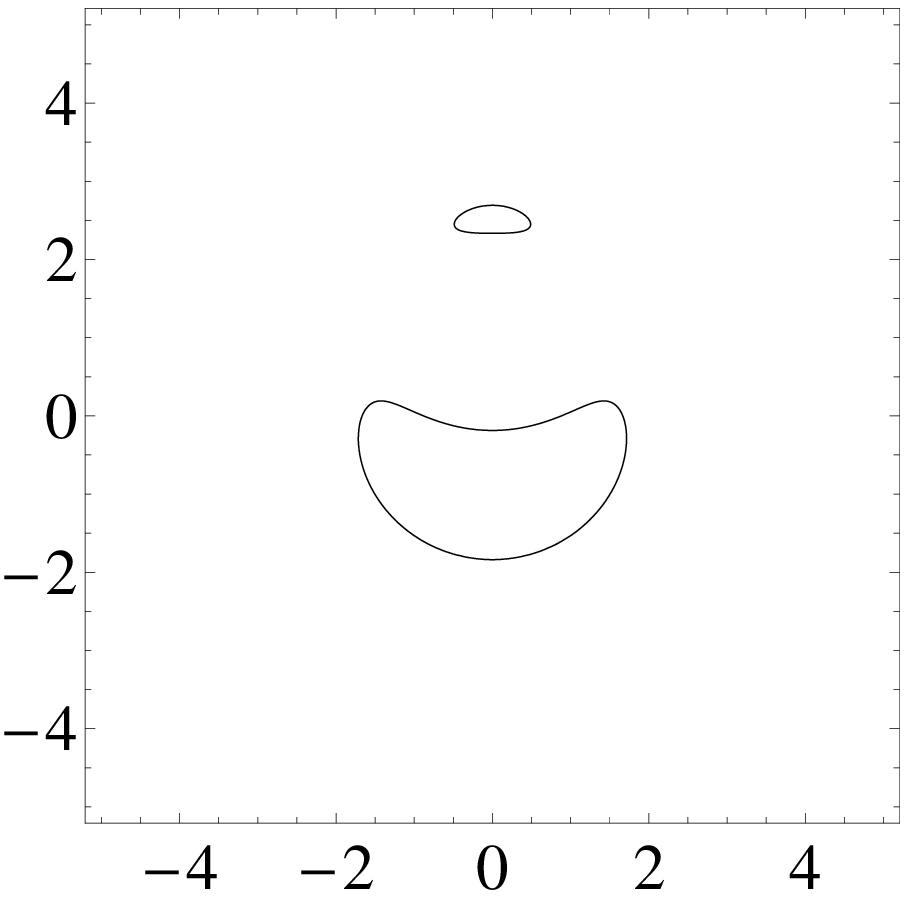} \\
            $t=1/4$;\  &
            $t=1/3$;\  &
            $t=2/3$ \\
        \end{tabular}}
\caption{\footnotesize{The curve $Re(H_{RSTUVW})=0$ in various ``time moments'' $t$ for $R(-1,0,0)$, $S(0,1,0)$, $T(1,0,0)$, $U(-2,0,i)$, $V(0,-1,-i)$ and $W(2,0,i)$. } }
        \label{GTUVW11}
\end{figure}

\begin{figure}[h]
        \setlength{\tabcolsep}{ 0 pt }{\scriptsize\tt
        \begin{tabular}{ ccc }
            \includegraphics[width=0.33\textwidth]{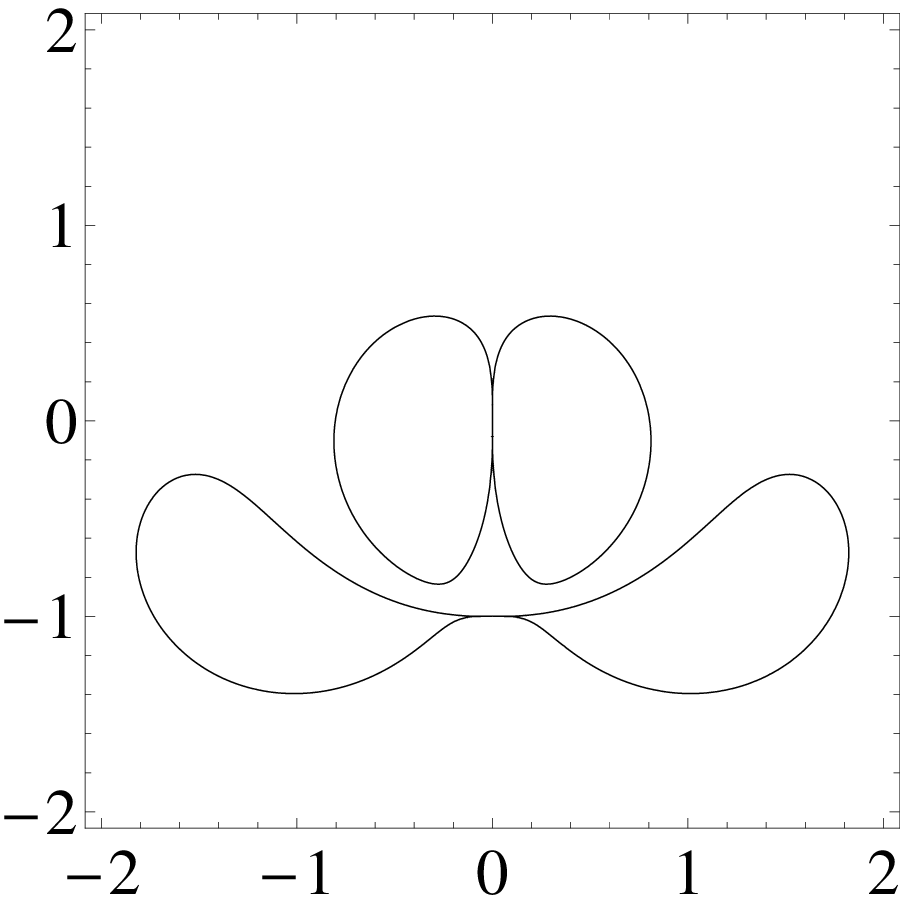} &
            \includegraphics[width=0.33\textwidth]{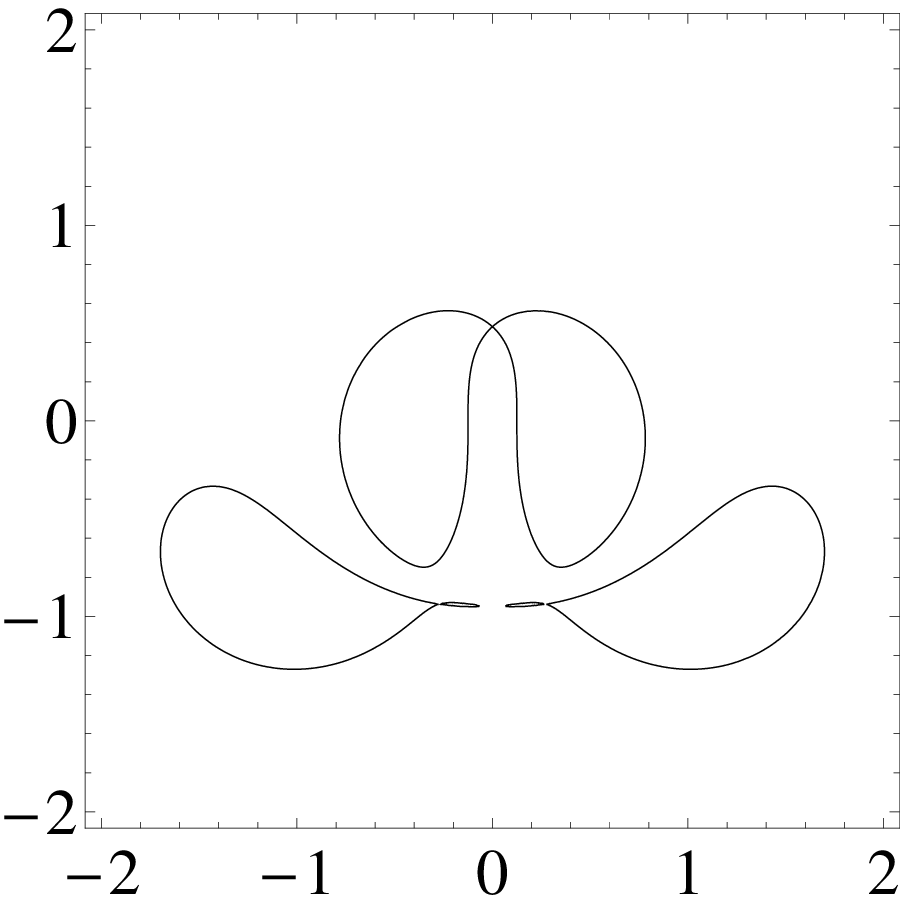} &
            \includegraphics[width=0.33\textwidth]{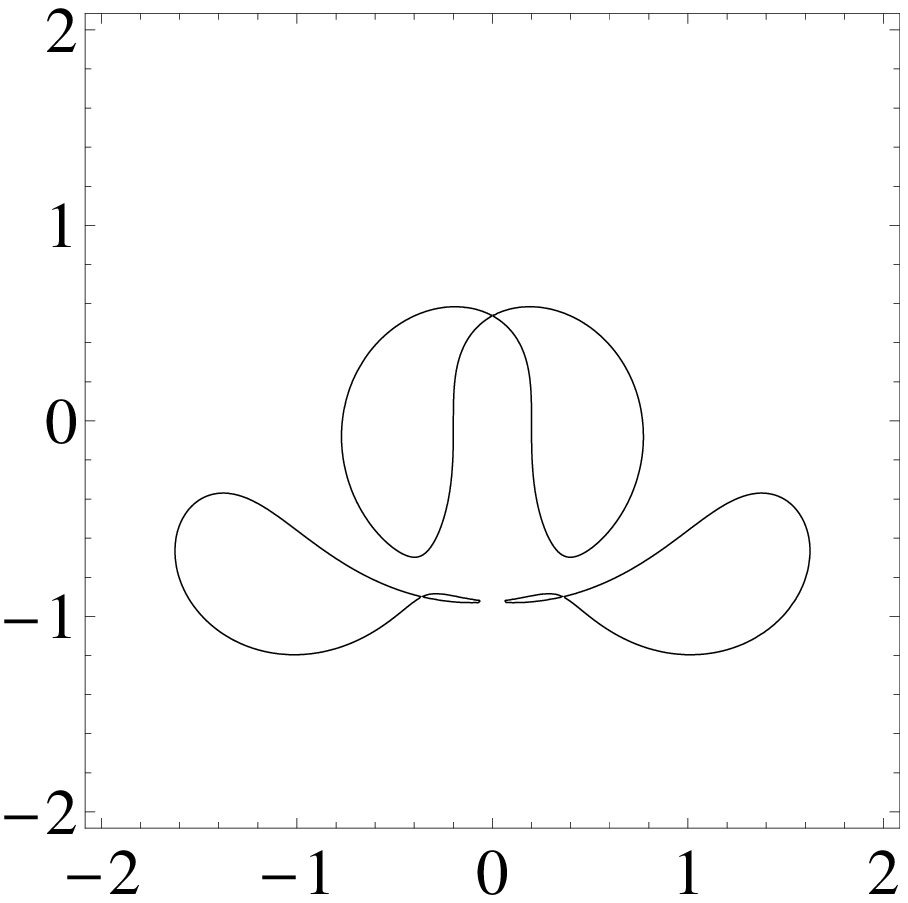} \\
            $t=-1$ &
            $t=-7/8$ &
            $t=-4/5$ \\
            \includegraphics[width=0.33\textwidth]{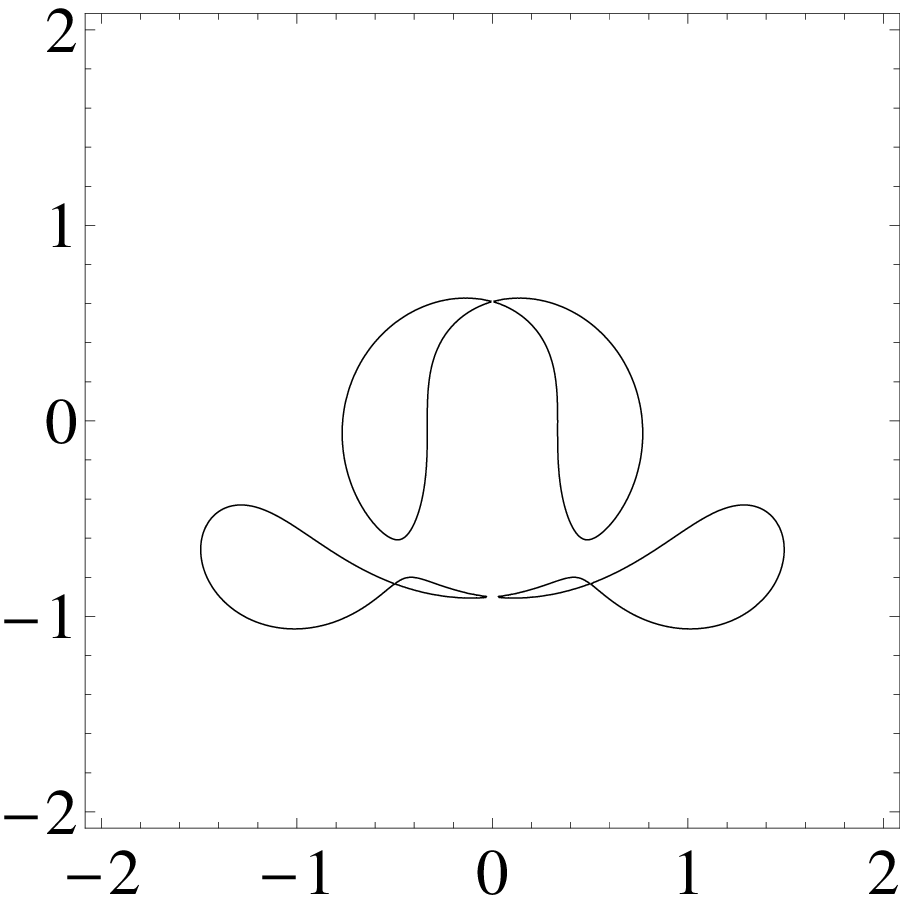} &
            \includegraphics[width=0.33\textwidth]{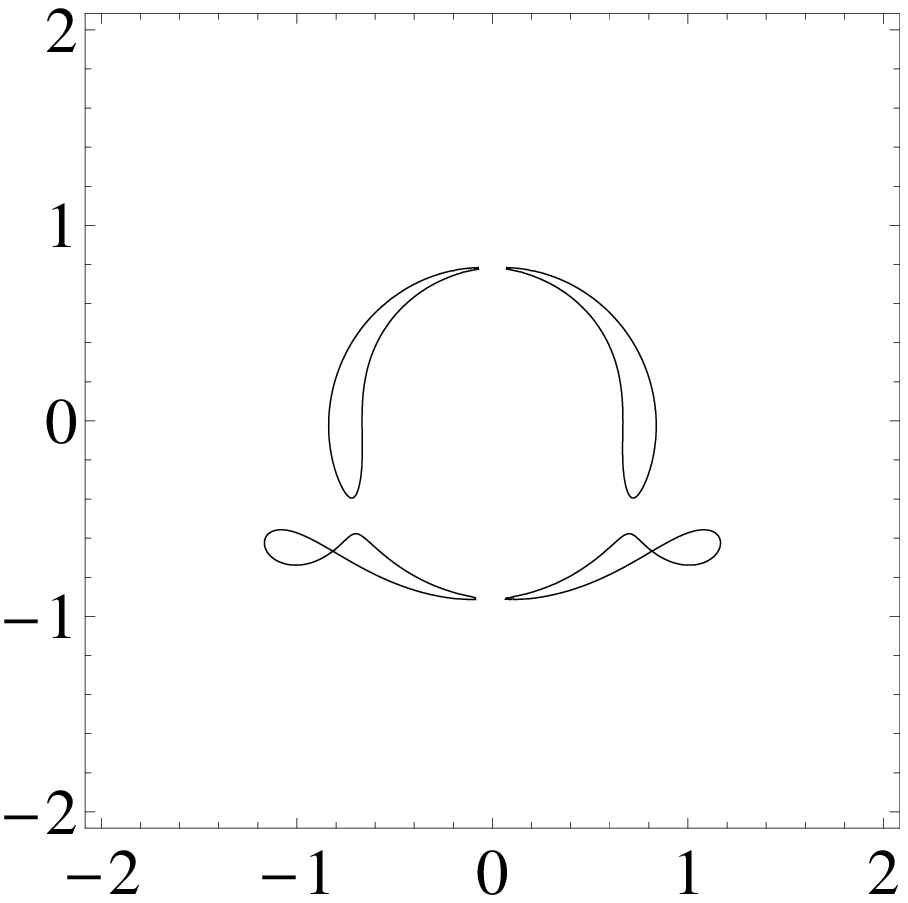} &
            \includegraphics[width=0.33\textwidth]{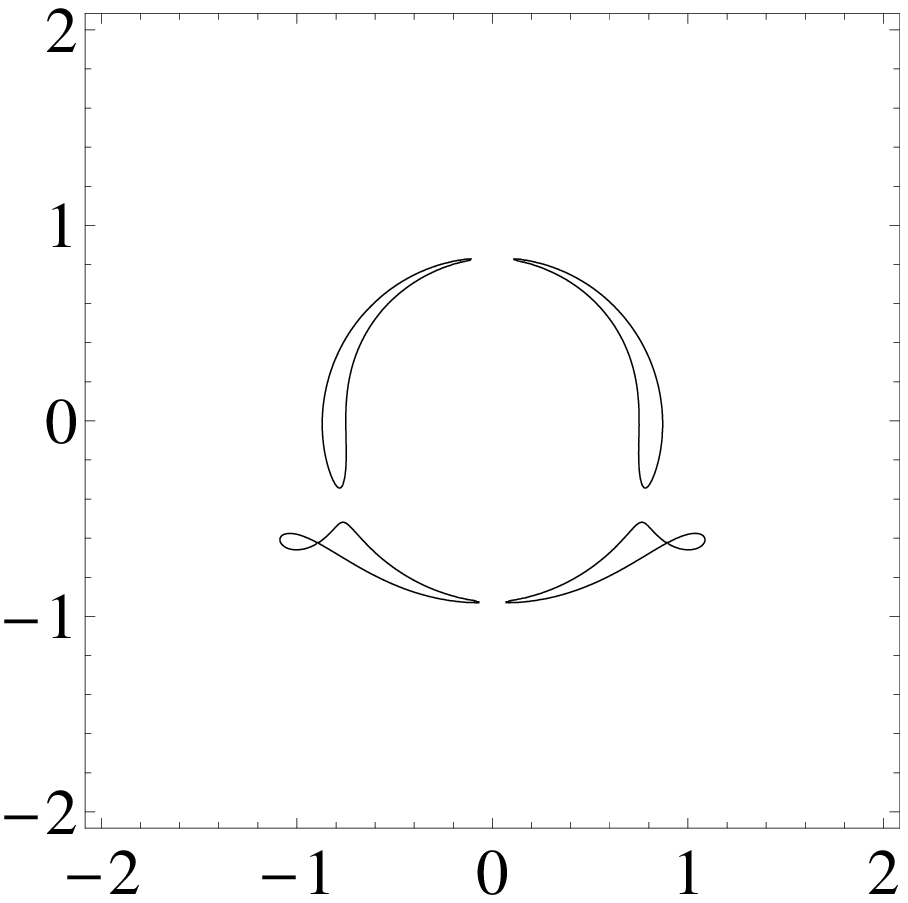} \\
            $t=-2/3$;\  &
            $t=-1/3$;\  &
            $t=-1/4$ \\
            \includegraphics[width=0.33\textwidth]{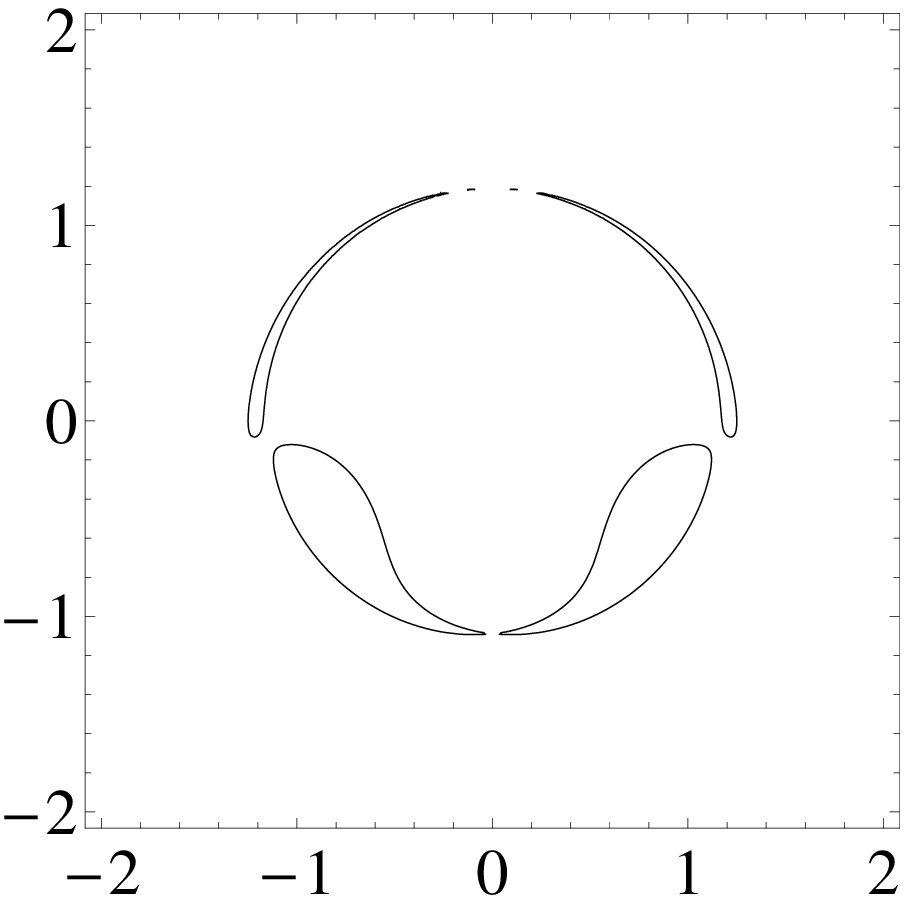} &
            \includegraphics[width=0.33\textwidth]{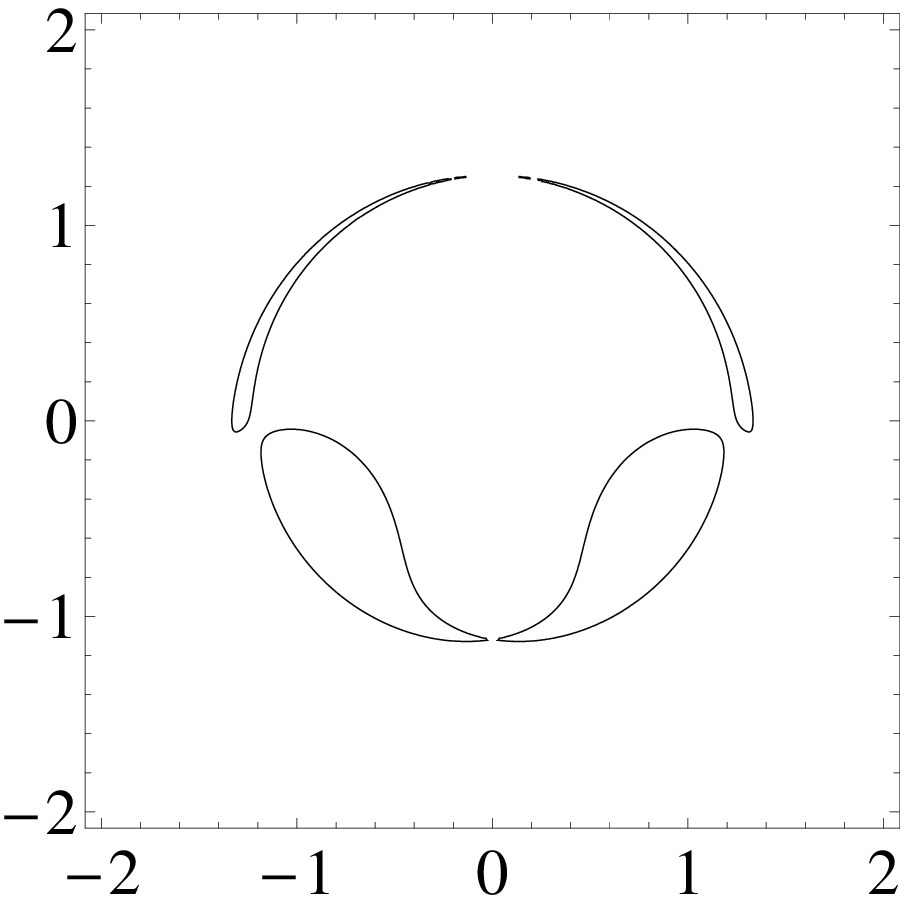} &
            \includegraphics[width=0.33\textwidth]{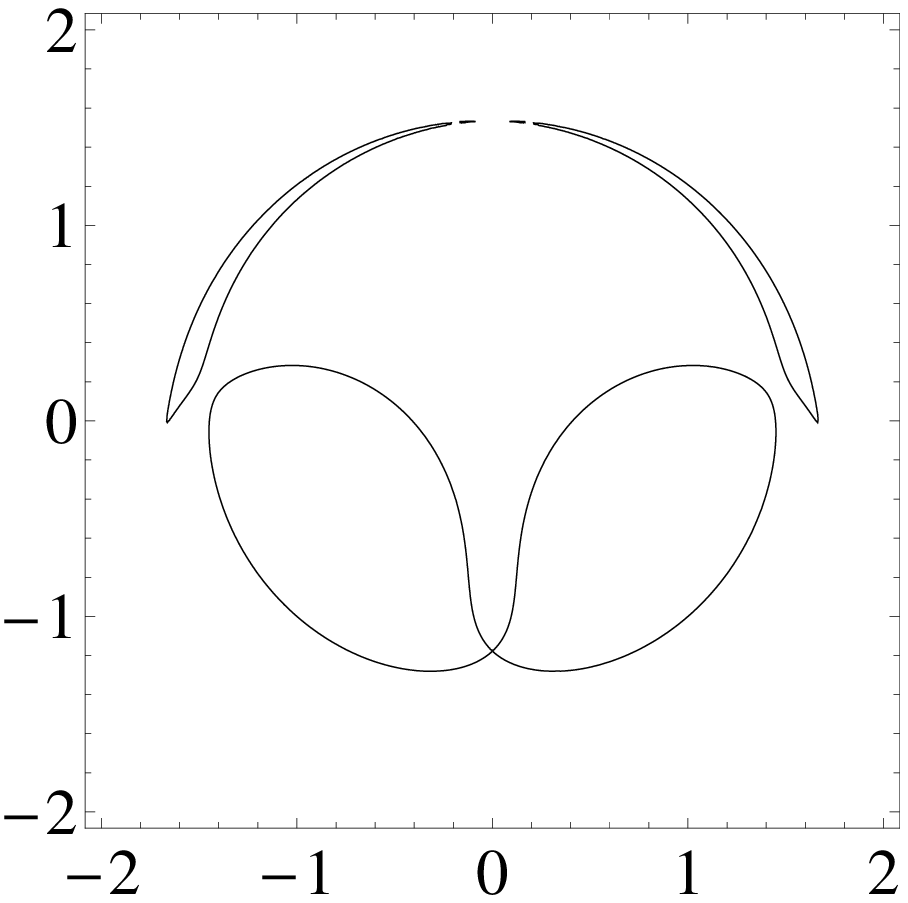} \\
            $t=1/4$;\  &
            $t=1/3$;\  &
            $t=2/3$ \\
        \end{tabular}}
\caption{\footnotesize{The curve $Re(G_{RSTUVW})=0$ in various ``time moments'' $t$ for $R(-1,0,0)$, $S(0,1,0)$, $T(1,0,0)$, $U(-2,0,i)$, $V(0,-1,-i)$ and $W(2,0,i)$. } }
        \label{GTUVW12}
\end{figure}

\begin{figure}[h]
        \setlength{\tabcolsep}{ 0 pt }{\scriptsize\tt
        \begin{tabular}{ ccc }
            \includegraphics[width=0.33\textwidth]{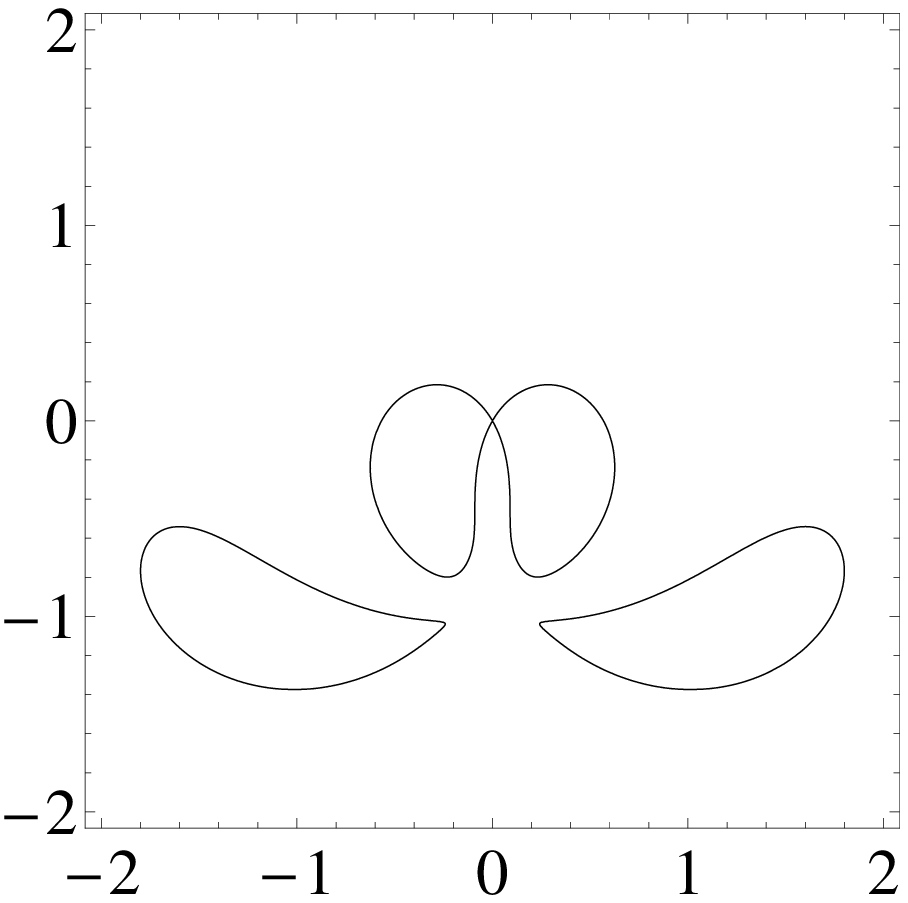} &
            \includegraphics[width=0.33\textwidth]{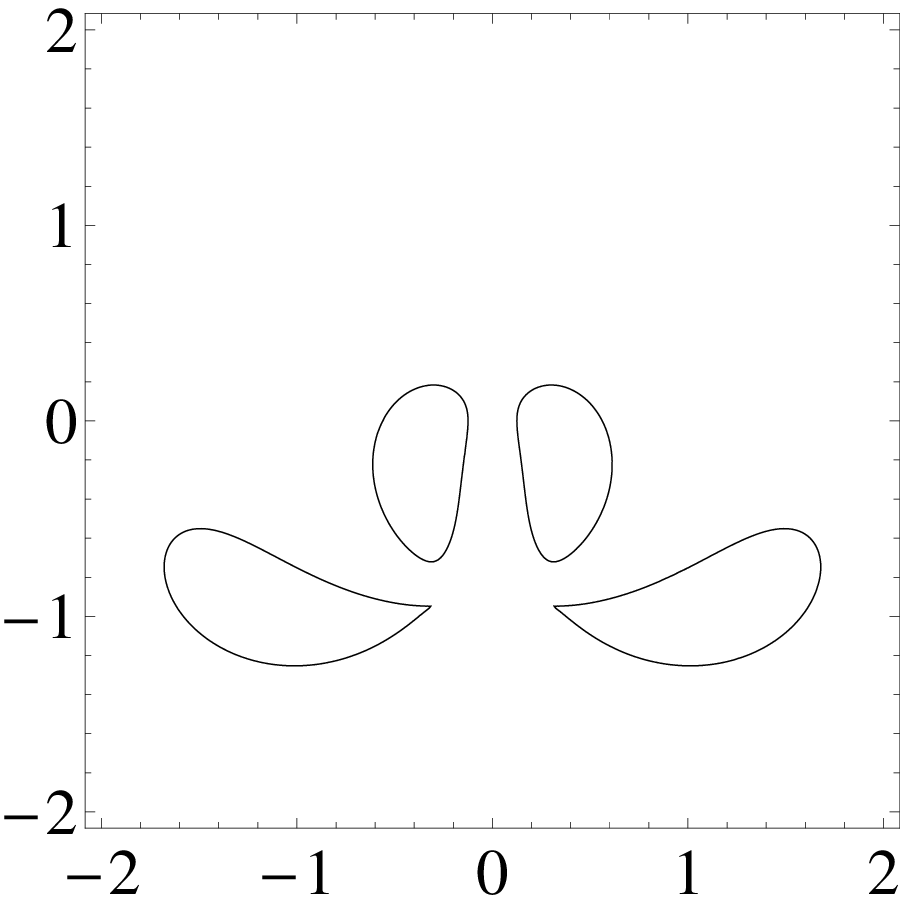} &
            \includegraphics[width=0.33\textwidth]{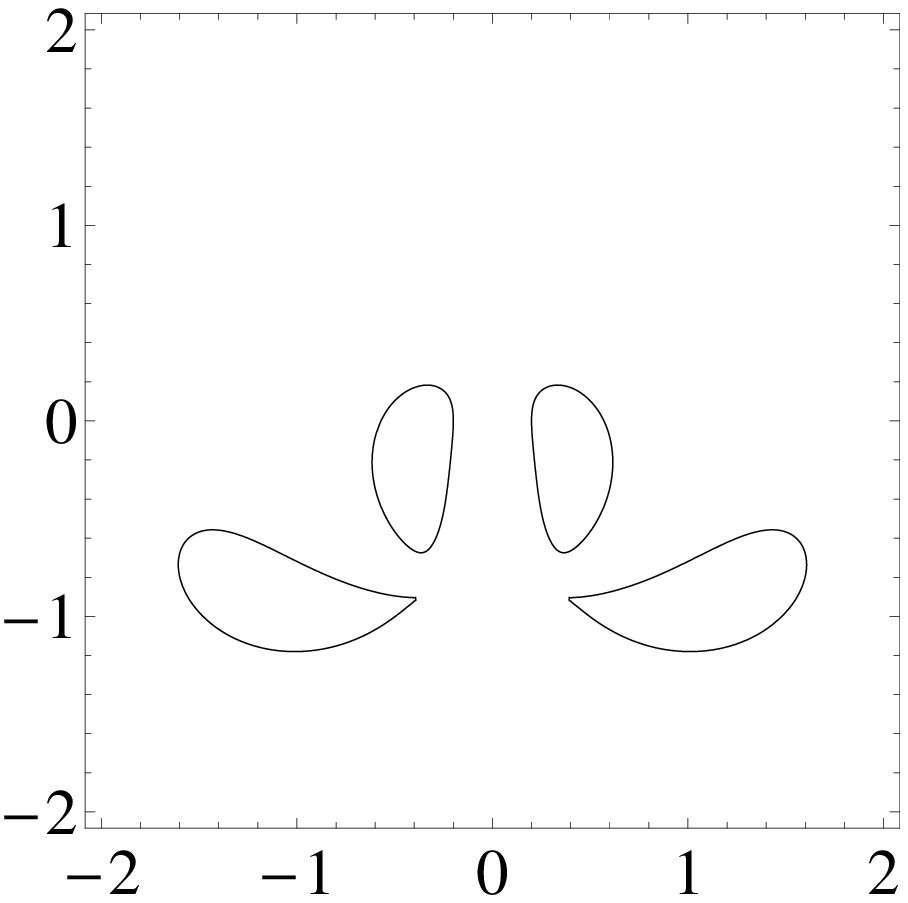} \\
            $t=-1$ &
            $t=-7/8$ &
            $t=-4/5$ \\
            \includegraphics[width=0.33\textwidth]{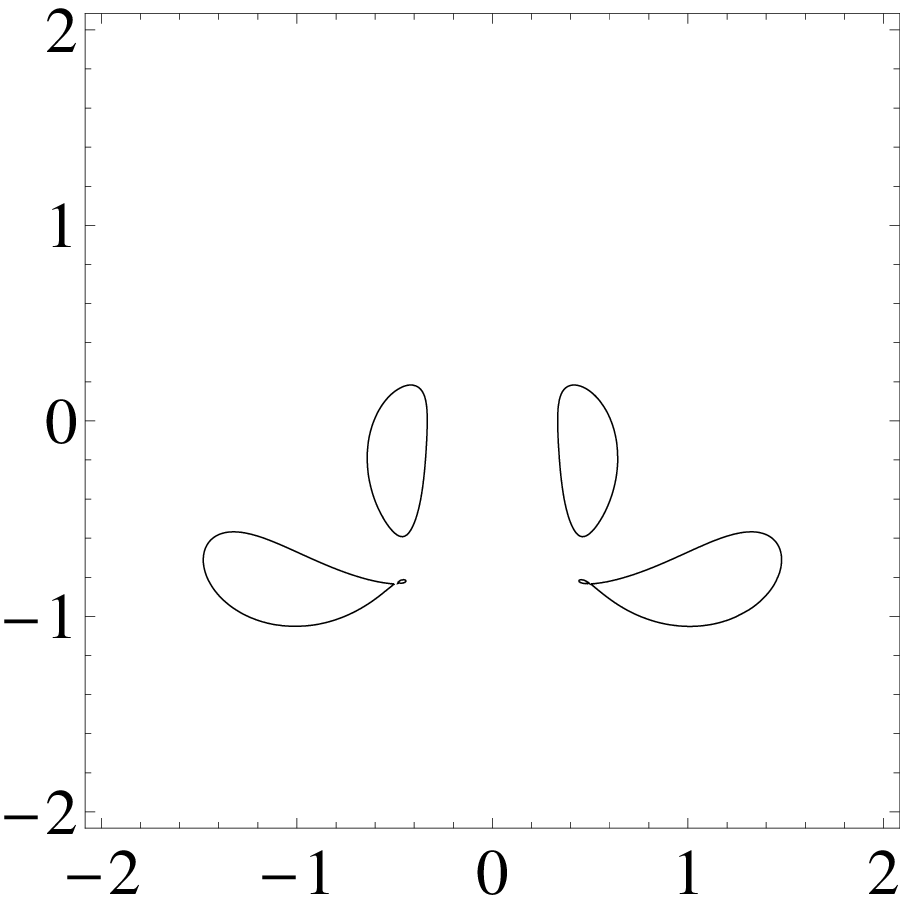} &
            \includegraphics[width=0.33\textwidth]{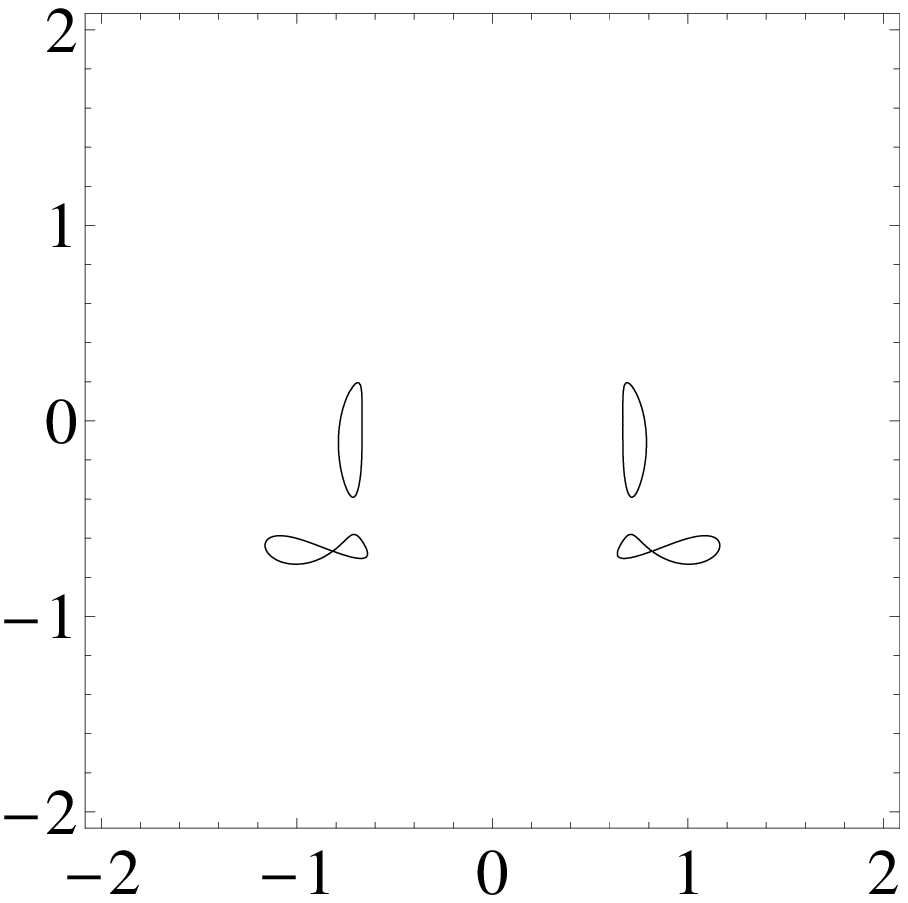} &
            \includegraphics[width=0.33\textwidth]{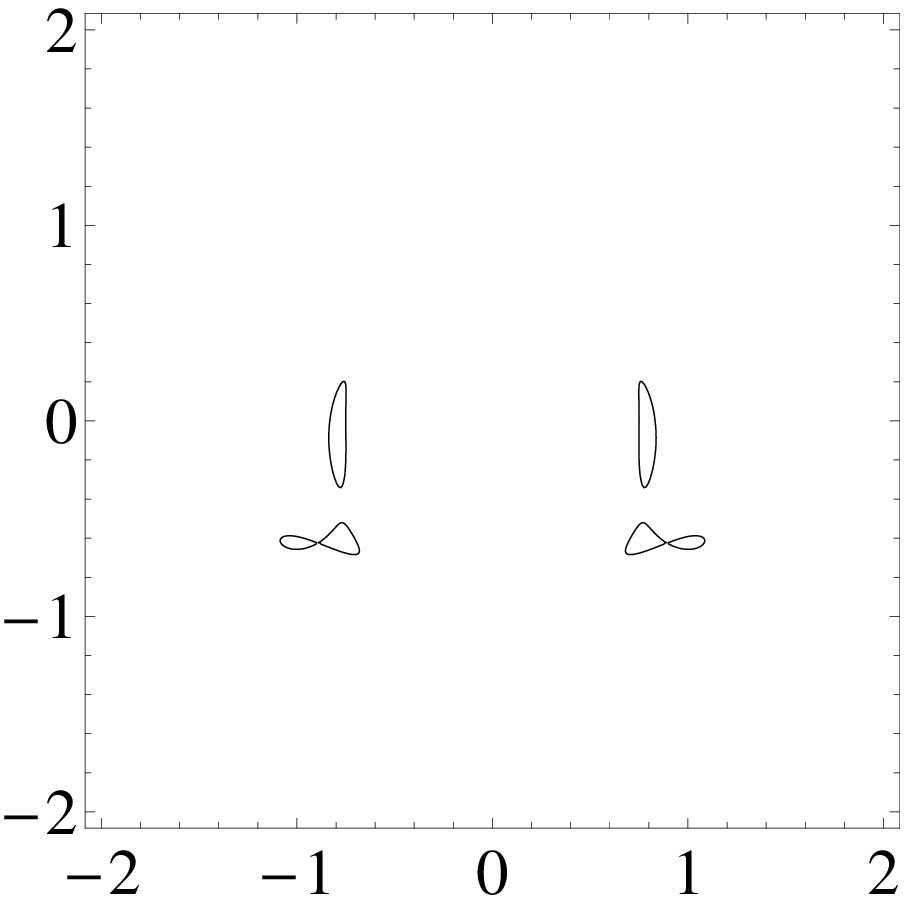} \\
            $t=-2/3$;\  &
            $t=-1/3$;\  &
            $t=-1/4$ \\
            \includegraphics[width=0.33\textwidth]{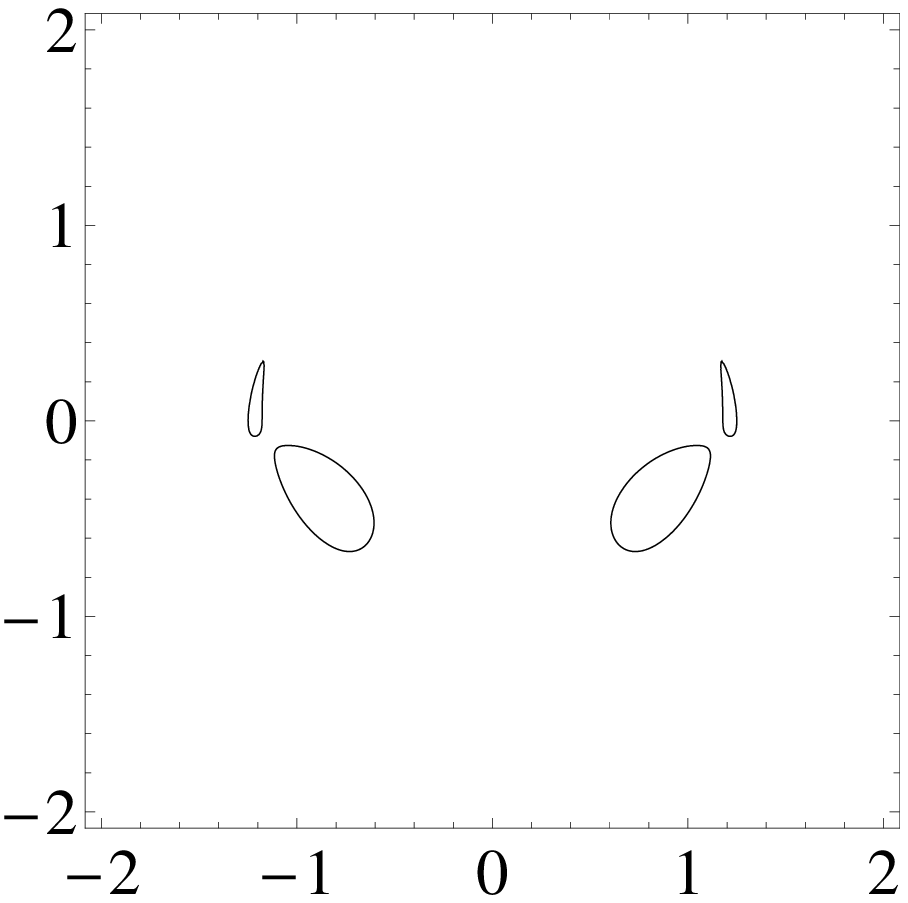} &
            \includegraphics[width=0.33\textwidth]{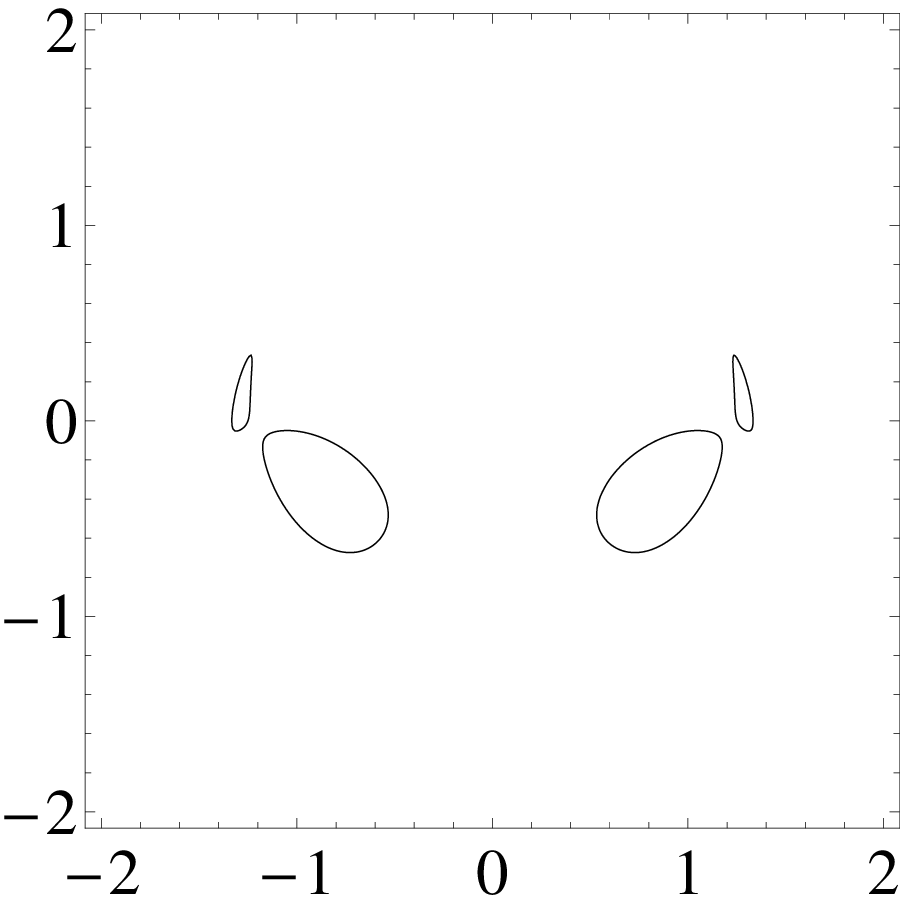} &
            \includegraphics[width=0.33\textwidth]{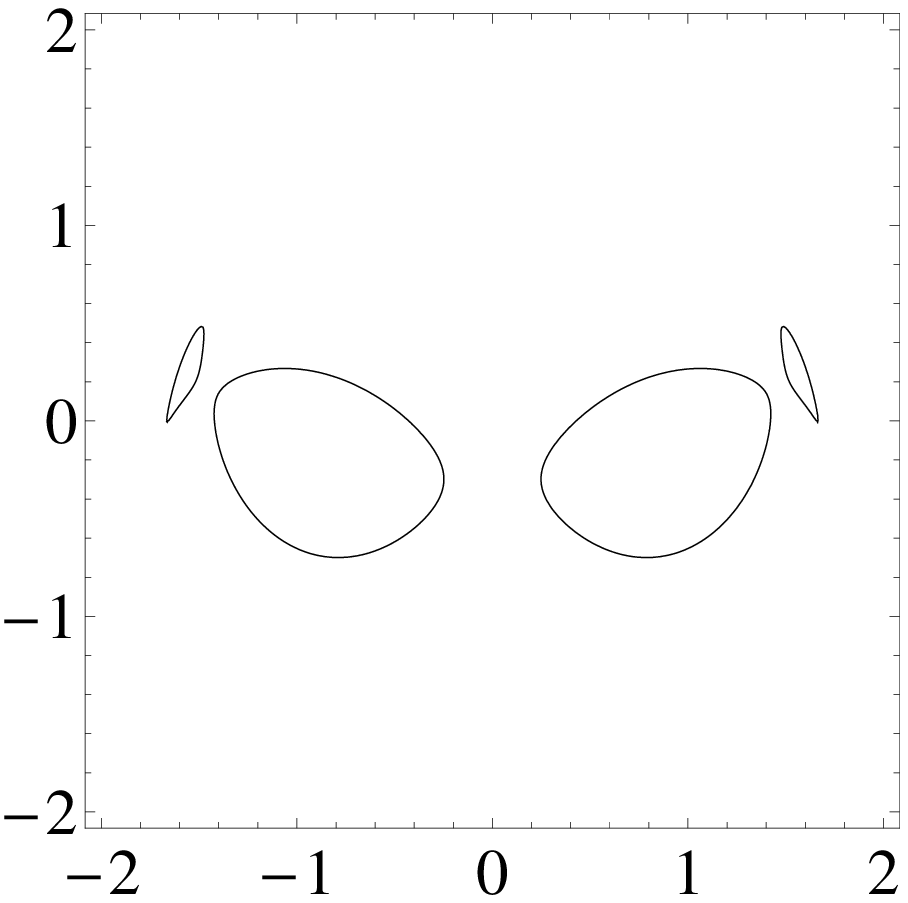} \\
            $t=1/4$;\  &
            $t=1/3$;\  &
            $t=2/3$ \\
        \end{tabular}}
\caption{\footnotesize{The curve $Re(F_{RSTUVW})=0$ for $\alpha=\pi/6$ rad in various ``time moments'' $t$ for $R(-1,0,0)$, $S(0,1,0)$, $T(1,0,0)$, $U(-2,0,i)$, $V(0,-1,-i)$ and $W(2,0,i)$. } }
        \label{GTUVW13}
\end{figure}

It is know that in the flat space with $N+1$ dimensions we can place as many $N$-spheres and $N$-simplexes as we wish. It is also known that each such $N$-sphere is defined uniquely by means of $N+2$ of its points, and each $N$-simplex in the flat space with $N+1$ dimensions has $N+1$ vertexes \cite{Sommerville}-\cite{Berger}. It is clear that a family of $N$-spheres passes through the vertexes of one such $N$-simplex. If we fix two distinct $N$-simplexes in our $(N+1)$-dimensional flat space, then two distinct families of $N$-spheres will pass through their vertexes. The elements of these families can have one common tangent point, or intersect in $(N-1)$-dimensional sphere, or have no common points. In the present paper we shall consider the couples of $N$-spheres of two distinct families of $N$-spheres, passing through the vertexes of two given $N$-simplexes, that have at least one common point. Moreover, for the angle $\alpha$ between the vectors connecting the common points of the couples of the $N$-spheres from the first and the second family with their centers, we assume that $\cos^{2}\alpha=const$. In the case when the two $N$-spheres touch at a single point it lies on the line connecting their centers, and then $\alpha$ is either $0$ or $\pi$ rad. Our goal is, given the value of $\cos^{2}\alpha$, to find the geometric place of the common points of the defined set of couples intersecting $N$-spheres in the Euclidean or pseudo-Euclidean space with $N+1$ dimensions.\\
For simplicity, we shall first consider the one-dimensional case. There are given two families of circles in the flat $2$-dimensional space, in which we have a fixed Cartesian coordinate system with arbitrarily fixed center $O$. The first family of circles passes through the points $T(x_{T},y_{T})$ и $U(x_{U},y_{U})$ (the vertexes of the first one-dimensional simplex $TU$), and the second one through the points $V(x_{V},y_{V})$ и $W(x_{W},y_{W})$ (the vertexes of the second one-dimensional simplex $VW$). We choose one common point of a given circle from the first family and the corresponding circle from the second family. Let us denote this common point by $K(x_{K},y_{K})$ (with radius vector $\vec r_{K}=(x_{K},y_{K})$). Then the equations of these two circles are \cite{Beyer}
\begin{equation}\label{circlesKTUVVW}
  \left| \begin{array}{cccc} x^{2}+y^{2} & x & y & 1 \\ x_{K}^{2}+y_{K}^{2} & x_{K} & y_{K} & 1 \\ x_{T}^{2}+y_{T}^{2} & x_{T} & y_{T} & 1 \\ x_{U}^{2}+y_{U}^{2} & x_{U} & y_{U} & 1 \\ \end{array} \right|=0,\
  \left| \begin{array}{cccc} x^{2}+y^{2} & x & y & 1 \\ x_{K}^{2}+y_{K}^{2} & x_{K} & y_{K} & 1 \\ x_{V}^{2}+y_{V}^{2} & x_{V} & y_{V} & 1 \\ x_{W}^{2}+y_{W}^{2} & x_{W} & y_{W} & 1 \\ \end{array} \right|=0.
\end{equation}
The equation of the first circle can also be written in the form
\begin{equation}\label{circleKTU0}
  \left(x-\frac{A_{KTU}}{2S_{KTU}}\right)^{2}+\left(y-\frac{B_{KTU}}{2S_{KTU}}\right)^{2}=
  \frac{1}{S_{KTU}}\left(\frac{A_{KTU}^{2}+B_{KTU}^{2}}{4S_{KTU}}+C_{KTU}\right),
\end{equation}
where
\begin{equation}\label{SCKTU}
  S_{KTU}\equiv \left| \begin{array}{ccc}x_{K} & y_{K} & 1 \\ x_{T} & y_{T} & 1 \\ x_{U} & y_{U} & 1 \\ \end{array} \right|,\
  C_{KTU}\equiv \left| \begin{array}{ccc}x_{K}^{2}+y_{K}^{2} & x_{K} & y_{K} \\ x_{T}^{2}+y_{T} ^{2}& x_{T} & y_{T} \\ x_{U}^{2}+y_{U}^{2} & x_{U} & y_{U} \\ \end{array} \right|,
\end{equation}
\begin{equation}\label{BAKTU}
  B_{KTU}\equiv -\left| \begin{array}{ccc}x_{K}^{2}+y_{K}^{2} & x_{K} & 1 \\ x_{T}^{2}+y_{T} ^{2}& x_{T} & 1 \\ x_{U}^{2}+y_{U}^{2} & x_{U} & 1 \\ \end{array} \right|,\
  A_{KTU}\equiv \left| \begin{array}{ccc}x_{K}^{2}+y_{K}^{2} & y_{K} & 1 \\ x_{T}^{2}+y_{T} ^{2}& y_{T} & 1 \\ x_{U}^{2}+y_{U}^{2} & y_{U} & 1 \\ \end{array} \right|.
\end{equation}
Equation (\ref{circleKTU0}) implies that the coordinates of the center of the circle are
\begin{equation}\label{x0KTU}
  {}_{0}x_{KTU}=\frac{A_{KTU}}{2S_{KTU}}=\frac{1}{2}\left| \begin{array}{ccc}x_{K}^{2}+y_{K}^{2} & y_{K} & 1 \\ x_{T}^{2}+y_{T} ^{2}& y_{T} & 1 \\ x_{U}^{2}+y_{U}^{2} & y_{U} & 1 \\ \end{array} \right|\left| \begin{array}{ccc}x_{K} & y_{K} & 1 \\ x_{T} & y_{T} & 1 \\ x_{U} & y_{U} & 1 \\ \end{array} \right|^{-1},
\end{equation}
\begin{equation}\label{y0KTU}
  {}_{0}y_{KTU}=\frac{B_{KTU}}{2S_{KTU}}=-\frac{1}{2}\left| \begin{array}{ccc}x_{K}^{2}+y_{K}^{2} & x_{K} & 1 \\ x_{T}^{2}+y_{T} ^{2}& x_{T} & 1 \\ x_{U}^{2}+y_{U}^{2} & x_{U} & 1 \\ \end{array} \right|\left| \begin{array}{ccc}x_{K} & y_{K} & 1 \\ x_{T} & y_{T} & 1 \\ x_{U} & y_{U} & 1 \\ \end{array} \right|^{-1}.
\end{equation}
For center of the other circle we have similar relations:
\begin{equation}\label{x0KVW}
  {}_{0}x_{KVW}=\frac{A_{KVW}}{2S_{KVW}}=\frac{1}{2}\left| \begin{array}{ccc}x_{K}^{2}+y_{K}^{2} & y_{K} & 1 \\ x_{V}^{2}+y_{V} ^{2}& y_{V} & 1 \\ x_{W}^{2}+y_{W}^{2} & y_{W} & 1 \\ \end{array} \right|\left| \begin{array}{ccc}x_{K} & y_{K} & 1 \\ x_{V} & y_{V} & 1 \\ x_{W} & y_{W} & 1 \\ \end{array} \right|^{-1},
\end{equation}
\begin{equation}\label{y0KVW}
  {}_{0}y_{KVW}=\frac{B_{KVW}}{2S_{KVW}}=-\frac{1}{2}\left| \begin{array}{ccc}x_{K}^{2}+y_{K}^{2} & x_{K} & 1 \\ x_{V}^{2}+y_{V} ^{2}& x_{V} & 1 \\ x_{W}^{2}+y_{W}^{2} & x_{W} & 1 \\ \end{array} \right|\left| \begin{array}{ccc}x_{K} & y_{K} & 1 \\ x_{V} & y_{V} & 1 \\ x_{W} & y_{W} & 1 \\ \end{array} \right|^{-1}.
\end{equation}
The vectors $\vec r_{K}-{}_{0}\vec r_{KTU}$ and $\vec r_{K}-{}_{0}\vec r_{KVW}$ enclose a fixed angle $\alpha$. Therefore for the common points of the considered two families of circles, we can write the following equation
\begin{equation}\label{eqvecrK}
  (\cos\alpha)|\vec r_{K}-{}_{0}\vec r_{KTU}||\vec r_{K}-{}_{0}\vec r_{KVW}|-(\vec r_{K}-{}_{0}\vec r_{KTU})\cdot(\vec r_{K}-{}_{0}\vec r_{KVW})=0.
\end{equation}
If we consider in detail the differences between these radius vectors, we shall see that
\begin{equation}\label{rKminusr0KTU}
  \vec r_{K}-{}_{0}\vec r_{KTU}=\frac{1}{2S_{KTU}}\ (a_{TU}(x_{K},y_{K}),b_{TU}(x_{K},y_{K})),
\end{equation}
\begin{equation}\label{rKminusr0KVW}
  \vec r_{K}-{}_{0}\vec r_{KVW}=\frac{1}{2S_{KVW}}\ (a_{VW}(x_{K},y_{K}),b_{VW}(x_{K},y_{K})),
\end{equation}
where
\begin{equation}\label{abCD}
  a_{CD}(x,y)=\left| \begin{array}{cccc} 2x & 1 & 0 & 0 \\ x^{2}+y^{2} & x & y & 1 \\ x_{C}^{2}+y_{C}^{2} & x_{C} & y_{C} & 1 \\ x_{D}^{2}+y_{D}^{2} & x_{D} & y_{D} & 1 \\ \end{array} \right|,\
  b_{CD}(x,y)=\left| \begin{array}{cccc} 2y & 0 & 1 & 0 \\ x^{2}+y^{2} & x & y & 1 \\ x_{C}^{2}+y_{C}^{2} & x_{C} & y_{C} & 1 \\ x_{D}^{2}+y_{D}^{2} & x_{D} & y_{D} & 1 \\ \end{array} \right|.
\end{equation}
Thus (\ref{eqvecrK}) is reduced to
\begin{equation}\label{eqrK2}
  (\cos^{2}\alpha)\left(a^{2}_{TU}+b^{2}_{TU}\right)\left(a^{2}_{VW}+b^{2}_{VW}\right)
  -\left(a_{TU}a_{VW}+b_{TU}b_{VW}\right)^{2}=0.
\end{equation}
It is convenient to utilize the identity
\begin{equation}\label{XY2}
  \left(X_{1}^{2}+X_{2}^{2}\right)\left(Y_{1}^{2}+Y_{2}^{2}\right)-\left(X_{1}Y_{1}+X_{2}Y_{2}\right)^{2}=
  \left| \begin{array}{cc} X_{1} & X_{2} \\ Y_{1} & Y_{2} \end{array}\right|^{2}.
\end{equation}
Through it (\ref{eqrK2}) takes the form
\begin{equation}\label{eqrK2a}
  F_{TUVW}(\alpha,x,y)\equiv\left| \begin{array}{cc} a_{TU} & b_{TU} \\ a_{VW} & b_{VW} \end{array}\right|^{2}-
  (\sin^{2}\alpha)\left(a^{2}_{TU}+b^{2}_{TU}\right)\left(a^{2}_{VW}+b^{2}_{VW}\right)=0.
\end{equation}
In the case of tangent circles we have that the considered geometric place of their tangent points is given by the equation $F_{TUVW}(k\pi,x,y)=0$ ($k=0,1$), which is now of the form
\begin{equation}\label{eqrK2tangential}
  G_{TUVW}(x,y)\equiv\left| \begin{array}{cc} a_{TU} & b_{TU} \\ a_{VW} & b_{VW} \end{array}\right|=0.
\end{equation}
In view of (\ref{abCD}), the function $G_{TUVW}(x,y)$ can be written explicitly in the form
\begin{equation}\label{funcGTUVW}
  \begin{array}{c} G_{TUVW}(x,y)\equiv\left| \begin{array}{cccc} 2x & 1 & 0 & 0 \\ x^{2}+y^{2} & x & y & 1 \\ x_{T}^{2}+y_{T}^{2} & x_{T} & y_{T} & 1 \\ x_{U}^{2}+y_{U}^{2} & x_{U} & y_{U} & 1 \\ \end{array} \right|\left| \begin{array}{cccc} 2y & 0 & 1 & 0 \\ x^{2}+y^{2} & x & y & 1 \\ x_{V}^{2}+y_{V}^{2} & x_{V} & y_{V} & 1 \\ x_{W}^{2}+y_{W}^{2} & x_{W} & y_{W} & 1 \\ \end{array} \right|-\\ \\ \ \ \ \ \ \ \ \ \ \ \ \ \ \ \ \ \ -\left| \begin{array}{cccc} 2y & 0 & 1 & 0 \\ x^{2}+y^{2} & x & y & 1 \\ x_{T}^{2}+y_{T}^{2} & x_{T} & y_{T} & 1 \\ x_{U}^{2}+y_{U}^{2} & x_{U} & y_{U} & 1 \\ \end{array} \right|\left| \begin{array}{cccc} 2x & 1 & 0 & 0 \\ x^{2}+y^{2} & x & y & 1 \\ x_{V}^{2}+y_{V}^{2} & x_{V} & y_{V} & 1 \\ x_{W}^{2}+y_{W}^{2} & x_{W} & y_{W} & 1 \\ \end{array} \right|. \end{array}
\end{equation}
We observe that the geometric place of the points of tangency of two families of circles, of which the members of the first family pass through the points $T(x_{T},y_{T})$ and $U(x_{U},y_{U})$, and the members of the second one through the points $V(x_{V},y_{V})$ and $W(x_{W},y_{W})$, is an algebraic curve of degree four. In the general case of two families of intersecting circles in the flat two-dimensional space, for which the angle $\alpha$ between the vectors, connecting the common points of the two circles in a couple with their centers, we have that the corresponding geometric place of the common points of these two families of circles is generally given by an algebraic equation of degree eight - (\ref{eqrK2a}). When $\alpha=\pi/2$ rad, this equation can again be represented as an algebraic equation of degree four as for $\alpha$ equal to zero or $\pi$ rad. Indeed, for $\alpha=\pi/2$ rad, (\ref{eqrK2}) implies
\begin{equation}\label{eqrK2orthogonal}
    H_{TUVW}(x,y)\equiv a_{TU}a_{VW}+b_{TU}b_{VW}=0.
\end{equation}
Let us explicitly write this function. We have
\begin{equation}\label{funcHTUVW}
  \begin{array}{c} H_{TUVW}(x,y)\equiv\left| \begin{array}{cccc} 2x & 1 & 0 & 0 \\ x^{2}+y^{2} & x & y & 1 \\ x_{T}^{2}+y_{T}^{2} & x_{T} & y_{T} & 1 \\ x_{U}^{2}+y_{U}^{2} & x_{U} & y_{U} & 1 \\ \end{array} \right|\left| \begin{array}{cccc} 2x & 1 & 0 & 0 \\ x^{2}+y^{2} & x & y & 1 \\ x_{V}^{2}+y_{V}^{2} & x_{V} & y_{V} & 1 \\ x_{W}^{2}+y_{W}^{2} & x_{W} & y_{W} & 1 \\ \end{array} \right|+\\ \\ \ \ \ \ \ \ \ \ \ \ \ \ \ \ \ \ \ +\left| \begin{array}{cccc} 2y & 0 & 1 & 0 \\ x^{2}+y^{2} & x & y & 1 \\ x_{T}^{2}+y_{T}^{2} & x_{T} & y_{T} & 1 \\ x_{U}^{2}+y_{U}^{2} & x_{U} & y_{U} & 1 \\ \end{array} \right|\left| \begin{array}{cccc} 2y & 0 & 1 & 0 \\ x^{2}+y^{2} & x & y & 1 \\ x_{V}^{2}+y_{V}^{2} & x_{V} & y_{V} & 1 \\ x_{W}^{2}+y_{W}^{2} & x_{W} & y_{W} & 1 \\ \end{array} \right|. \end{array}
\end{equation}
Thus the problem we consider is solved in the one-dimensional case. The curves given by equation (\ref{eqrK2tangential}) in the cases when $TU\nparallel VW$ and we cannot circumscribe a circle around $TUVW$ are closed and have no points of self-intersection. When $TU\parallel VW$ and we cannot circumscribe a circle around $TUVW$, the curves are open and they are asymptotically close to lines at infinity. The case when we can circumscribe a circle around $TUVW$ is degenerate---the curves decompose into two circles, one of which is circumscribed around $TUVW$ provided that it is not a rectangle or a bilateral trapezium. Moreover, the two circles coincide when $TU$ and $VW$ are the diagonals of $TUVW$. Finally, if $TUVW$ is a rectangle or a bilateral trapezium, then the curves decompose into a line and the circle circumscribed around $TUVW$. Figures \ref{GTUVW1} and \ref{GTUVW2} show various non-degenerate curves (\ref{eqrK2tangential}). Figure \ref{GTUVW4} shows the degenerate case when a circle can be circumscribed around $TUVW$, and also a similar non-degenerate variant of curve (\ref{eqrK2tangential}). Figure \ref{GTUVW6} shows the degenerate case in which $TUVW$ is rectangle or a bilateral trapezium ($TU\parallel VW$). Figures \ref{GTUVW7a} and \ref{GTUVW7} show examples of curve (\ref{eqrK2orthogonal}), and Figure \ref{GTUVW8} contains illustrations of the curves defined by equation (\ref{eqrK2a}). Generally speaking, equation (\ref{eqrK2orthogonal}) defines a class of curves, whose graphs consist of two parts. Each ot these parts is closed except when the angle between $TU$ and $VW$ is $\pi/2$ rad. Curve (\ref{eqrK2orthogonal}) for $\angle(TU,VW)=\pi/2$ rad is shown on Figure \ref{GTUVW7a}. A peculiar feature of this case is that we have a degeneration of the curve when one of the segments $TU$ or $VW$ lies on the perpendicular bisector of the other. Then the curve degenerate into a circle, passing through two of the four points and the perpendicular bisector of these two points, provided that $TU$ and $VW$ are not the diagonals of a square or another quadrangle round which a circle can be circumscribed. In these cases we have degeneration again---curve (\ref{eqrK2orthogonal}) reduces to the lines through $TU$ and $VW$, or to two circles that intersect at a right angle and pass through $TU$ and $VW$.\\
The graphs of the curves defined by equation (\ref{eqrK2a}) have two, three (in one special case), or four parts, which, as in the previous case, are closed, provided that the angle between $TU$ and $VW$ is not equal to $\alpha$. An example of curve (\ref{eqrK2a}) for $\alpha=\angle(TU,VW)$ is given on Figure \ref{GTUVW8a}. The special case is realized when the segments $TU$ and $VW$ connect the midpoints of the opposite sides of a parallelogram. Then the graph of curve (\ref{eqrK2a}) has three parts for $\alpha=\angle(TU,VW)$. Figure \ref{GTUVW8b} shows an example of this peculiar case of curve (\ref{eqrK2a}).\\
When $N>1$ we can conduct similar considerations. The equations of the two families of $N$-spheres are stated by means of the condition the determinants of two square matrixes of order $N+3$ to be equal to zero,
\begin{equation}\label{spheresKSymplex1Symplex2}
    \begin{array}{ll}
    \left| \begin{array}{cccccc} \sum\limits_{i=1}^{n}x_{i}^{2} & x_{1} & x_{2} & ... & x_{n} & 1 \\ \sum\limits_{i=1}^{n}x_{iK}^{2} & x_{1K} & x_{2K} & ... & x_{nK} & 1 \\ \sum\limits_{i=1}^{n}x_{iS_{1}}^{2} & x_{1S_{1}} & x_{2S_{1}}& ... & x_{nS_{1}} & 1 \\ ... & ... & ... & ... & ... & 1 \\ \sum\limits_{i=1}^{n}x_{iS_{n}}^{2} & x_{1S_{n}} & x_{2S_{n}}& ... & x_{nS_{n}} & 1 \end{array} \right|=0,\\
    \\
    \left| \begin{array}{cccccc} \sum\limits_{i=1}^{n}x_{i}^{2} & x_{1} & x_{2} & ... & x_{n} & 1 \\ \sum\limits_{i=1}^{n}x_{iK}^{2} & x_{1K} & x_{2K} & ... & x_{nK} & 1 \\ \sum\limits_{i=1}^{n}x_{iV_{1}}^{2} & x_{1V_{1}} & x_{2V_{1}}& ... & x_{nV_{1}} & 1 \\ ... & ... & ... & ... & ... & 1 \\ \sum\limits_{i=1}^{n}x_{iV_{n}}^{2} & x_{1V_{n}} & x_{2V_{n}}& ... & x_{nV_{n}} & 1 \end{array} \right|=0.
    \end{array}
\end{equation}
Here $n=N+1$, and $x_{iK}$ ($i=1,2,...,n$) are the coordinates of any common point of the couples of $N$-spheres in the containing $N+1$-dimensional flat space. As in the one-dimensional case, we denote this point by $K$. Similarly, $x_{iS_{j}}$ and $x_{iV_{j}}$ ($i=1,2,...,n$; $j=1,2,...,n$) are the coordinates of the $j$-th vertex of the first and the second $N$-dimensional simplex in the same $N+1$-dimensional flat space. We denote the vertexes (as points) of these simplexes by $S_{j}$ and $V_{j}$. Let the radius vectors of the centers of any representatives of the couples of $N$-spheres be ${}_{0}\vec r_{KS_{1}..S_{n}}$ and ${}_{0}\vec r_{KV_{1}..V_{n}}$. About the differences between the radius vectors of the point $K$ and these two vectors, we can write relations similar to (\ref{rKminusr0KTU}), (\ref{rKminusr0KVW}) and (\ref{abCD}). We have
\begin{equation}\label{rKminusr0KSj}
  \vec r_{K}-{}_{0}\vec r_{KS_{1}...S_{n}}=\frac{1}{2V_{KS_{1}...S_{n}}}\ (a_{1S_{1}...S_{n}}(\vec r_{K}),a_{2S_{1}...S_{n}}(\vec r_{K}),...,a_{nS_{1}...S_{n}}(\vec r_{K})),
\end{equation}
\begin{equation}\label{rKminusr0KVj}
  \vec r_{K}-{}_{0}\vec r_{KV_{1}...V_{n}}=\frac{1}{2V_{KV_{1}...V_{n}}}\ (a_{1V_{1}...V_{n}}(\vec r_{K}),a_{2V_{1}...V_{n}}(\vec r_{K}),...,a_{nV_{1}...V_{n}}(\vec r_{K})),
\end{equation}
where
\begin{equation}\label{VKiSjVj}
    \begin{array}{ll}
    V_{KS_{1}...S_{n}}=\left| \begin{array}{ccccc} x_{1K} & x_{2K} & ... & x_{nK} & 1 \\ x_{1S_{1}} & x_{2S_{1}}& ... & x_{nS_{1}} & 1 \\ ... & ... & ... & ... & 1 \\ x_{1S_{n}} & x_{2S_{n}}& ... & x_{nS_{n}} & 1 \end{array} \right|,\\
    \\
     V_{KV_{1}...V_{n}}=\left| \begin{array}{ccccc} x_{1K} & x_{2K} & ... & x_{nK} & 1 \\ x_{1V_{1}} & x_{2V_{1}}& ... & x_{nV_{1}} & 1 \\ ... & ... & ... & ... & 1 \\ x_{1V_{n}} & x_{2V_{n}}& ... & x_{nV_{n}} & 1 \end{array} \right|,
    \end{array}
\end{equation}
and
\begin{equation}\label{aiBj}
    \begin{array}{ll}
    a_{1B_{1}...B_{n}}(\vec r)=\left| \begin{array}{cccccc} 2x_{1} & 1 & 0 & ... & 0 & 0 \\ \sum\limits_{i=1}^{n}x_{i}^{2} & x_{1} & x_{2} & ... & x_{n} & 1 \\ \sum\limits_{i=1}^{n}x_{iB_{1}}^{2} & x_{1B_{1}} & x_{2B_{1}}& ... & x_{nB_{1}} & 1 \\ ... & ... & ... & ... & ... & 1 \\ \sum\limits_{i=1}^{n}x_{iB_{n}}^{2} & x_{1B_{n}} & x_{2B_{n}}& ... & x_{nB_{n}} & 1 \end{array} \right|,\\
    \\
    a_{2B_{1}...B_{n}}(\vec r)=\left| \begin{array}{cccccc} 2x_{2} & 0 & 1 & ... & 0 & 0 \\ \sum\limits_{i=1}^{n}x_{i}^{2} & x_{1} & x_{2} & ... & x_{n} & 1 \\ \sum\limits_{i=1}^{n}x_{iB_{1}}^{2} & x_{1B_{1}} & x_{2B_{1}}& ... & x_{nB_{1}} & 1 \\ ... & ... & ... & ... & ... & 1 \\ \sum\limits_{i=1}^{n}x_{iB_{n}}^{2} & x_{1B_{n}} & x_{2B_{n}}& ... & x_{nB_{n}} & 1 \end{array} \right|,\\
    \\
    ...\ ,\\
    a_{nB_{1}...B_{n}}(\vec r)=\left| \begin{array}{cccccc} 2x_{n} & 0 & 0 & ... & 1 & 0 \\ \sum\limits_{i=1}^{n}x_{i}^{2} & x_{1} & x_{2} & ... & x_{n} & 1 \\ \sum\limits_{i=1}^{n}x_{iB_{1}}^{2} & x_{1B_{1}} & x_{2B_{1}}& ... & x_{nB_{1}} & 1 \\ ... & ... & ... & ... & ... & 1 \\ \sum\limits_{i=1}^{n}x_{iB_{n}}^{2} & x_{1B_{n}} & x_{2B_{n}}& ... & x_{nB_{n}} & 1 \end{array} \right|.
    \end{array}
\end{equation}
The equation, analogous to (\ref{eqrK2}), in the $N$-dimensional case is given by
\begin{equation}\label{eqrKN}
  (\cos^{2}\alpha)\left(\sum\limits_{i=1}^{n}a_{iS_{1}...S_{n}}^{2}\right)\left(\sum\limits_{i=1}^{n}a_{iV_{1}...V_{n}}^{2}\right)
  -\left(\sum\limits_{i=1}^{n}a_{iS_{1}...S_{n}}\ a_{iV_{1}...V_{n}}\right)^{2}=0.
\end{equation}
The identity we have to use now is
\begin{equation}\label{XYN}
  \left(\sum\limits_{i=1}^{n}X_{i}^{2}\right)\left(\sum\limits_{i=1}^{n}Y_{i}^{2}\right)
  -\left(\sum\limits_{i=1}^{n}X_{i}Y_{i}\right)^{2}=
  \sum\limits_{i=1}^{n-1}\sum\limits_{j>i}^{n}\left(\left| \begin{array}{cc} X_{i} & X_{j} \\ Y_{i} & Y_{j} \end{array}\right|^{2}\right).
\end{equation}
Thus (\ref{eqrKN}) is reduced to
\begin{equation}\label{eqrKNa}
    \begin{array}{ll}
    F_{S_{1}...S_{n}V_{1}...V_{n}}(\alpha, \vec r)\equiv\sum\limits_{i=1}^{n-1}\sum\limits_{j>i}^{n}\left(\left| \begin{array}{cc} a_{iS_{1}...S_{n}} & a_{jS_{1}...S_{n}} \\ a_{iV_{1}...V_{n}} & a_{jV_{1}...V_{n}} \end{array}\right|^{2}\right)-\\
    \\
    \ \ \ \ \ \ \ \ \ \ \ \ \ \ \ \ \ \ \ \ \ \ \ \ \ -(\sin^{2}\alpha)\left(\sum\limits_{i=1}^{n}a_{iS_{1}...S_{n}}^{2}\right)
    \left(\sum\limits_{i=1}^{n}a_{iV_{1}...V_{n}}^{2}\right)=0.
    \end{array}
\end{equation}
It is clear that unlike the one-dimensional case, the defined two families of $N$-spheres for $N>1$ do not determine an $N$-dimensional surface if $\alpha$ $0$ or $\pi$ rad because then there should hold
\begin{equation}\label{eqrKNtangential}
  G_{ijS_{1}...S_{n}V_{1}...V_{n}}(\vec r)\equiv\left| \begin{array}{cc} a_{iS_{1}...S_{n}} & a_{jS_{1}...S_{n}} \\ a_{iV_{1}...V_{n}} & a_{jV_{1}...V_{n}} \end{array}\right|=0,
\end{equation}
for every $i$ and $j$ ($j>i$). In the other cases this is not so. Moreover, just like in the one-dimensional case, the equation that describes the considered geometric place can be simplified for $\alpha=\pi/2$ rad. By (\ref{eqrKN}) we have
\begin{equation}\label{eqrKNorthogonal}
  H_{S_{1}...S_{n}V_{1}...V_{n}}(\vec r)\equiv\sum\limits_{i=1}^{n}a_{iS_{1}...S_{n}}\ a_{iV_{1}...V_{n}}=0.
\end{equation}
It is clear that for $n=2$ thie equation defines a surface in the flat three-dimensional space. Figure \ref{GTUVW9} shows an example of such a two-dimensional manifold. We see that it has the topology of a torus. Figure \ref{GTUVW10} shows an example of the general two-dimensional manifold (\ref{eqrKNa}) such that $\sin^{2}\alpha\neq 0$ or 1.\\
If we allow the containing flat $N+1$-dimensional space to be pseudo-Euclidean \cite{Greub}, \cite{NovikovFomenko}, then we can define the smooth manifolds, which are a multidimensional analog to the $G$-surfaces, given by (\ref{eqrK2tangential}) for $N=1$. Let $x_{nA}$ be $N+1$th coordinate of an arbitrary point $A$ of the containing flat space. In this case, in order that it is pseudo-Euclidean, it is sufficient to set
\begin{equation}\label{pseudoEuclidianCaseX}
  x_{nA}=it_{A};\ i^{2}\equiv -1,\ t_{A}\in\mathds{R}.
\end{equation}
Then $x_{N+1}=it$ and for $t=const\in\mathds{R}$ we can define the following two ``fixed-moment'' manifolds (along with their intersection):
\begin{equation}\label{eqrKNtangentialPseudoEuclidian}
    \begin{array}{ll}
    Re\left[\sum\limits_{k=1}^{n-1}\sum\limits_{l>k}^{n}G^{2}_{klS_{1}...S_{n}V_{1}...V_{n}}(t,\vec r)\right]=0,\\
    \\
    Im\left[\sum\limits_{k=1}^{n-1}\sum\limits_{l>k}^{n}G^{2}_{klS_{1}...S_{n}V_{1}...V_{n}}(t,\vec r)\right]=0.
    \end{array}
\end{equation}
Here $\vec r=\left(x_{1A},x_{2A},...,x_{NA}\right)\in\mathds{R}^{N}$. Of course, nothing prevents us from doing the same with (\ref{eqrKNorthogonal}) and (\ref{eqrKNa}). We have
\begin{equation}\label{eqrKNorthogonalPseudoEuclidian}
    \begin{array}{ll}
    Re\left[H_{S_{1}...S_{n}V_{1}...V_{n}}(t,\vec r)\right]=0,\\
    \\
    Im\left[H_{S_{1}...S_{n}V_{1}...V_{n}}(t,\vec r)\right]=0;
    \end{array}
\end{equation}
\begin{equation}\label{eqrKNaPseudoEuclidian}
    \begin{array}{ll}
    Re\left[F_{S_{1}...S_{n}V_{1}...V_{n}}(a,b,t,\vec r)\right]=0,\\
    \\
    Im\left[F_{S_{1}...S_{n}V_{1}...V_{n}}(a,b,t,\vec r)\right]=0,
    \end{array}
\end{equation}
where $\sin^{2}\alpha\equiv a+ib$ ($a=const\in\mathds{R},\ b=const\in\mathds{R}$). Figures \ref{GTUVW11} - \ref{GTUVW13} show examples of this type of manifolds. Equations (\ref{eqrKNtangentialPseudoEuclidian}), (\ref{eqrKNorthogonalPseudoEuclidian}) and (\ref{eqrKNaPseudoEuclidian}) implicitly give all possible smooth ``fixed-moment'' manifolds of the types we consider in the pseudo-Euclidean case with signature $\{n,1\}$ ($n=N+1$). The cases with different signature can be considered in a similar way.\\
All this manifolds are of interest, especially if we take into account that they possess quite noteworthy properties with regard to inversion.

\clearpage

\section*{Acknowledgments}
I would like to extend my gratitude to engineer Vladimir Borissovich Dryzhak. In the course of solving the problem he posed about the construction of a convex quadrangle in the plane on given feet of the altitudes from the intersection of its diagonals, I came to the generalization leading to the manifolds considered in the present work. I am thankful to Ivan Simeonov, who encouraged me during the preparation of the manuscript, and also to Kiril Delev, Prof. Kiril Bankov and Associate Prof. Borislav Draganov. Kiril Delev read the part about the $G$-curves, Prof. Bankov read the paper in its final form, and Prof. Draganov help me with English translation.

\end{document}